\documentclass[journal]{IEEEtran}
\usepackage{lipsum}
\usepackage{amsfonts}
\usepackage{graphicx}
\usepackage{epstopdf}
\usepackage{algorithmic}
\usepackage{subcaption}
\usepackage{amsmath}
\usepackage{hyperref}
\ifpdf
  \DeclareGraphicsExtensions{.eps,.pdf,.png,.jpg}
\else
  \DeclareGraphicsExtensions{.eps}
\fi

\usepackage{enumitem}
\setlist[enumerate]{leftmargin=.5in}
\setlist[itemize]{leftmargin=.5in}

\newtheorem{theorem}{Theorem}[section]

\newtheorem{lemma}[theorem]{Lemma}
\newtheorem{corollary}[theorem]{Corollary}

\newtheorem{definition}[theorem]{Definition}
\newtheorem{remark}[theorem]{Remark}

\ifCLASSINFOpdf
\else
\fi
\begin{document}
%
\title{Fast and Asymptotically Powerful Detection \\for Filamentary Objects in Digital Images}
%
%
%

\author{Kai~Ni,
        Shanshan~Cao,
        and~Xiaoming~Huo
\thanks{Kai Ni was with the School of Mathematics, Georgia Institute of Technology, Atlanta,
GA, 30332 USA, e-mail: kni0219@gmail.com.}
\thanks{Shanshan Cao is with Eli Lilly and Company, Indianapolis,
IN, 46225 USA, e-mail: shancao36@gmail.com, and this work was mostly done while she was in the Department of Industrial and System Engineering, Georgia Institute of Technology.}
\thanks{Xiaoming Huo is with the Department of Industrial and System Engineering, Georgia Institute of Technology, Atlanta,
GA, 30332 USA, e-mail: huo@gatech.edu.}
}

\maketitle

\begin{abstract}
Given an inhomogeneous chain embedded in a noisy image, we consider the conditions under which such an embedded chain is detectable.
Many applications, such as detecting moving objects, detecting ship wakes, can be abstracted as the detection on the existence of chains.
In this work, we provide the detection algorithm with low order of computation complexity  to detect the chain and the optimal theoretical detectability regarding SNR (signal to noise ratio) under the normal distribution model.
Specifically, we derive an analytical threshold that specifies what is detectable.
We design a longest significant chain detection algorithm, with computation complexity in the order of $O(n\log n)$.
We also prove that our proposed algorithm is asymptotically powerful, which means, as the dimension $n \rightarrow \infty$, the probability of false detection vanishes.
We further provide some simulated examples and a real data example, which validate our theory.
\end{abstract}


%
\IEEEpeerreviewmaketitle

\section{Introduction}
%
%
%
%
\IEEEPARstart{D}{etectability}  is a fundamental problem in many image processing tasks.
It is to determine whether detecting an object via a computer is doable.
Furthermore, if it is doable, what is an appropriate order of complexity for the associated algorithm.
In \cite{Donoho2005}, the authors proved a range of powerful results regarding the detection of the presence of a geometric object in an image with additive Gaussian noises; the efficient detection algorithm can have orders of complexity such as $O(n)$ or $O(n^2)$ depending on the class of the objects.
In \cite{Xiaoming2009}, a detection problem of unknown convex sets is presented; the infeasibility of adopting the generalized likelihood ratio test (which succeeded in \cite{Donoho2005}) is proved due to the large cardinality of the convex sets under consideration.
An approach via hv-parallelograms is analyzed by studying the minimax proportion of a hv-parallelogram inscribed in a convex set. By adopting such an approach, the authors present the efficiency and the lower order of complexity of the corresponding method.
In \cite{tandra2008snr}, the authors consider a broader question, which is called the ``SNR walls'', below which a detector will fail to be robust, no matter how long it can observe the channel, using the simple mathematical models for the uncertainty in the noise and fading processes.

Constantly advanced imaging technology and better software and hardware lead to demands and wishes to use digital images as a tool for evaluation and analysis.
In many applications, data or images collected by standard sensors (such as cameras and radars) are analyzed for the detection and recognition of targets.
Detecting an inhomogeneous region \cite{sutour2015estimation}, which can be modeled as a chain in a noisy environment \cite{Aos06Filament, Donoho2005}, is one of these problems. In applications of image detection problems, one class of questions is to determine whether or not some filamentary structures are present in the noisy picture.
There is a plethora of available statistical methods that can, in principle, be used for filaments detection and estimation.
These include: Principle curves in \cite{Hastie}, \cite{Kegl}, \cite{Sathyakama} and \cite{Smola}; nonparametric, penalized, maximum likelihood in \cite{RTibshirani}; parametric models in \cite{Stoica}; manifold learning techniques in \cite{Roweis}, \cite{JoshuaB} and \cite{HuoandChen}; gradient based methods in \cite{DNovikov} and \cite{Genovese}; methods from computational geometry in \cite{Dey}, \cite{Lee1999} and \cite{Cheng2005}; faint line segment detection in \cite{ExpandVision}; Ship Wakes ``V'' shape detection against a highly cluttered background in \cite {Vdetection} and underlying curvilinear structure in \cite{Aos06Filament}, \cite{AAOMT} and \cite{Donoho2005}.
See also \cite{perim1}, \cite{perim2} and \cite{perim3} for the applications of the percolation theory in this area.
Recently cluster detection problems have been fully studied in regular lattice (\cite{trailInMaze, ClustDetecPerco, AnomaClusterInNetwork}) and \cite{NetworksPoly} presented a similar model as ours.

In this paper, we consider the detection of chains with good continuation, which is defined later.
The problem of detecting a path in a network appears to represent a fundamental abstraction as many modern statistical detection problems can reasonably be formulated in this way.
In \cite{trailInMaze}, the authors consider a filamentary detection problem in a tree based graph with a very similar structure as ours.
In \cite{Aos06Filament} and \cite{NetworksPoly}, the authors use methods based on this structure model to detect graphs of  H$\ddot{o}$lder functions possibly hidden in point clouds.
In \cite{horev2015detection}, the authors develop the sublinear algorithm to detect long straight edges in large and noisy images by  processing only a small subset of the image pixels smartly.
In addition, they theoretically analyze the inevitable tradeoff between its detection performance and the allowed computational budget.
The essence of these problems is to detect whether a sequence of connected nodes (which are variables that can be measured) exhibit a peculiar behavior.
In \cite{WaterQuality}, for instance, the authors assess the water quality in a network of streams by performing a chemical analysis at various locations along the streams.
As a result, some locations are marked as problematic.
One may consider the set of all tested locations along the streams as nodes and  pairs of adjacent nodes located on the same stream, are connected, which yields the problem of detection of chains with good continuation.
One can then imagine that in order to trace the existence of a polluter, or look for the existence of an anomalous path along all streams, the problem would be to detect an upstream of a certain sensitive location.

\begin{figure}[htbp]
    \centering
    \begin{subfigure}[b]{0.48\textwidth}
        \includegraphics[width=\textwidth]{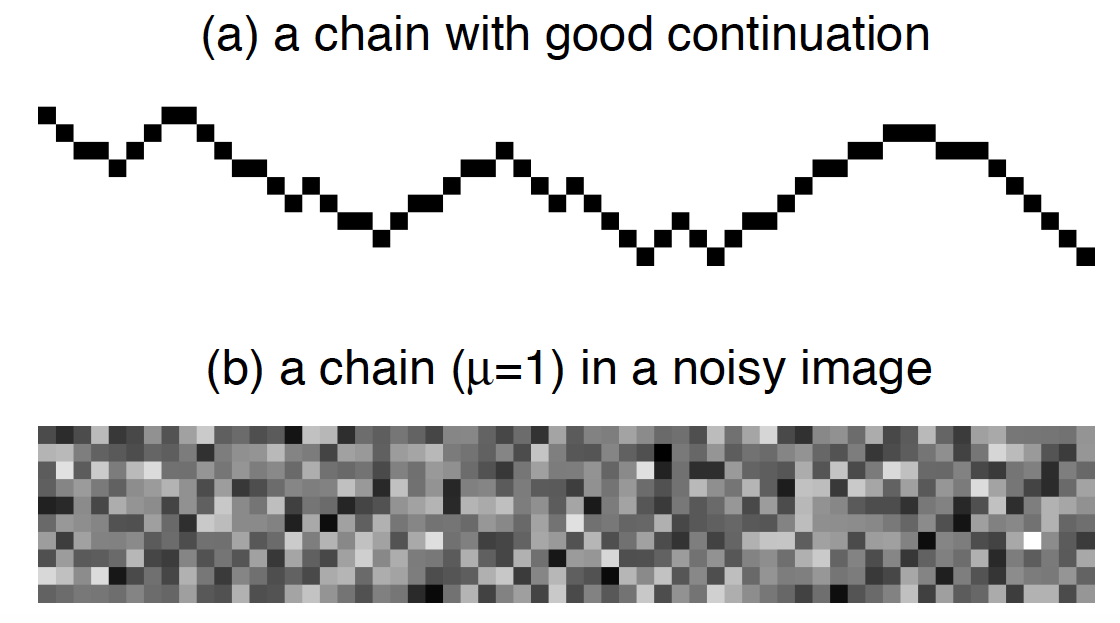}
    \end{subfigure}
    ~ 
    \begin{subfigure}[b]{0.48\textwidth}
        \includegraphics[width=\textwidth]{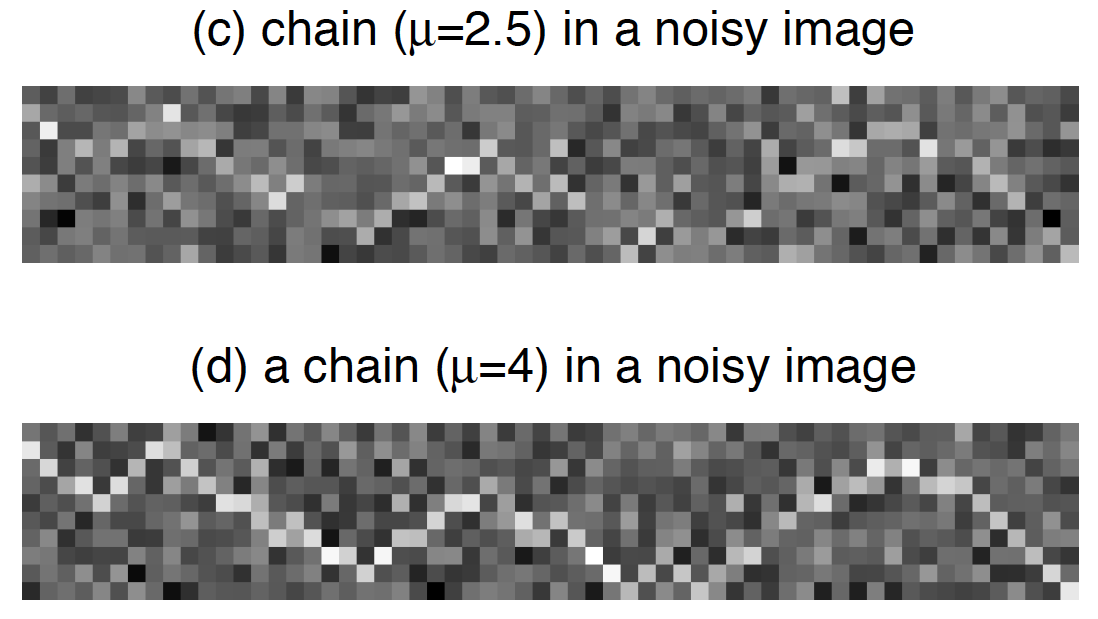}
    \end{subfigure}
    \caption{A chain with good continuation (a) and chains embedded in noisy images with means being equal to $1$ in (b), $2.5$ in (c), and $4.0$ in (d).}\label{compareFig}
\end{figure}

It is helpful to consider an illustration of our model at this point.
Fig. \ref{compareFig} (a) contains a chain with good continuation in an image with $10$-by-$60$ pixels; (b)-(d) present the same chain in noisy Gaussian random fields with different elevated means.
The detectability problem is to ask: when is the chain detectable and what is the order of complexity of the detection algorithm?
We have the following statistical formulation: the intensity at each pixel follows a normal distribution.
Inside the chain, the normal mean $\mu$ is a positive constant such as in Fig. \ref{compareFig} (b) $\mu=1.0$, (c) $\mu=2.5$ and (d) $\mu=4.0$, and the standard deviation $\sigma = 1$ is used in all the simulated examples in this work; while outside the chain, the normal means are $0$.
In (b) and (c), the chain can hardly be observed by eyes while in (d), the chain is clearly visible.

\begin{figure}[htbp]
    \centering
    \begin{subfigure}[b]{0.39\textwidth}
        \includegraphics[width=\textwidth]{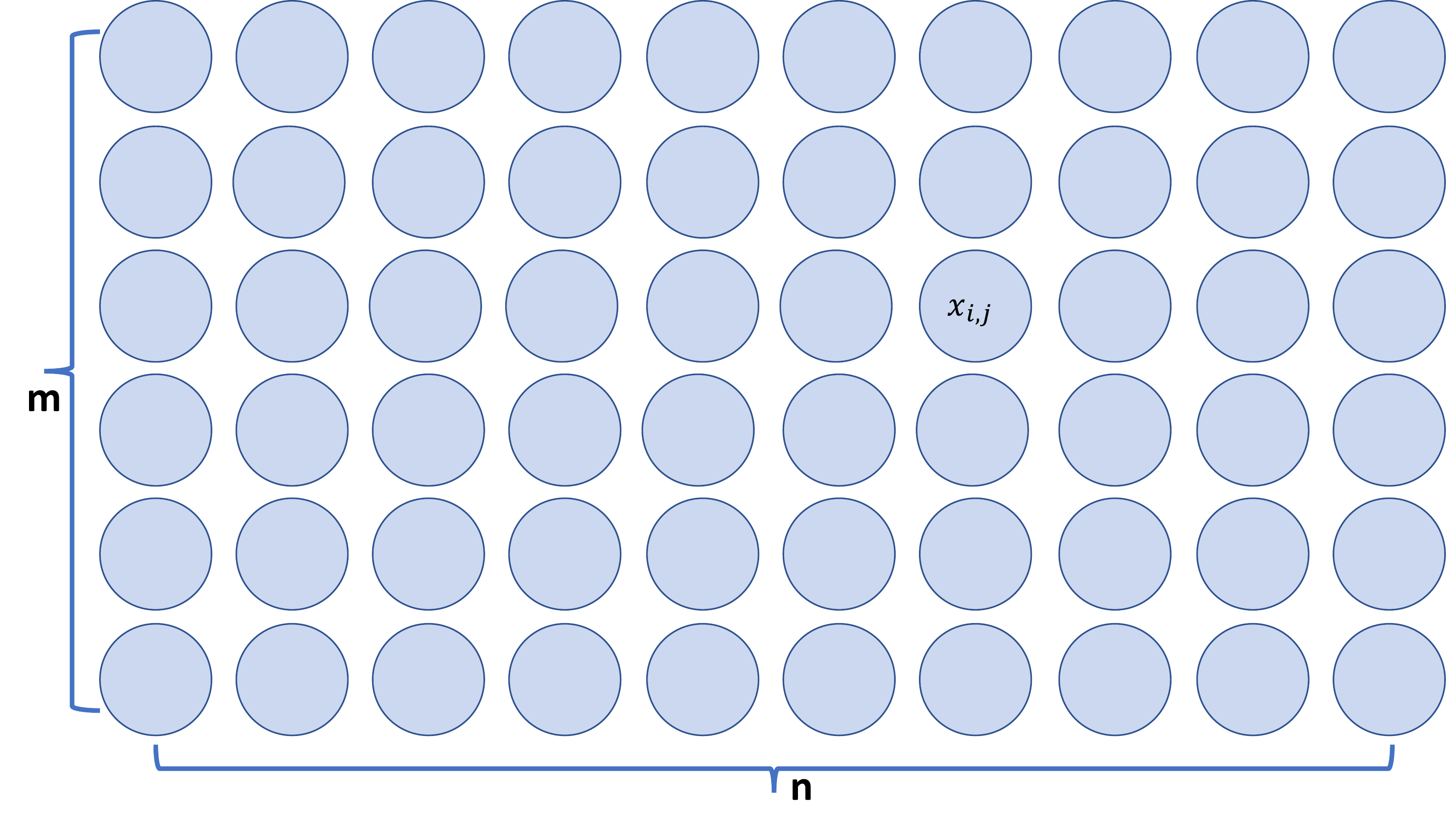}
\caption{The graph structure for an image, where the value of $x_{i,j}$ indicates the intensity of the corresponding node.}\label{fig:graph1}
    \end{subfigure}
    ~ 
    \begin{subfigure}[b]{0.39\textwidth}
        \includegraphics[width=\textwidth]{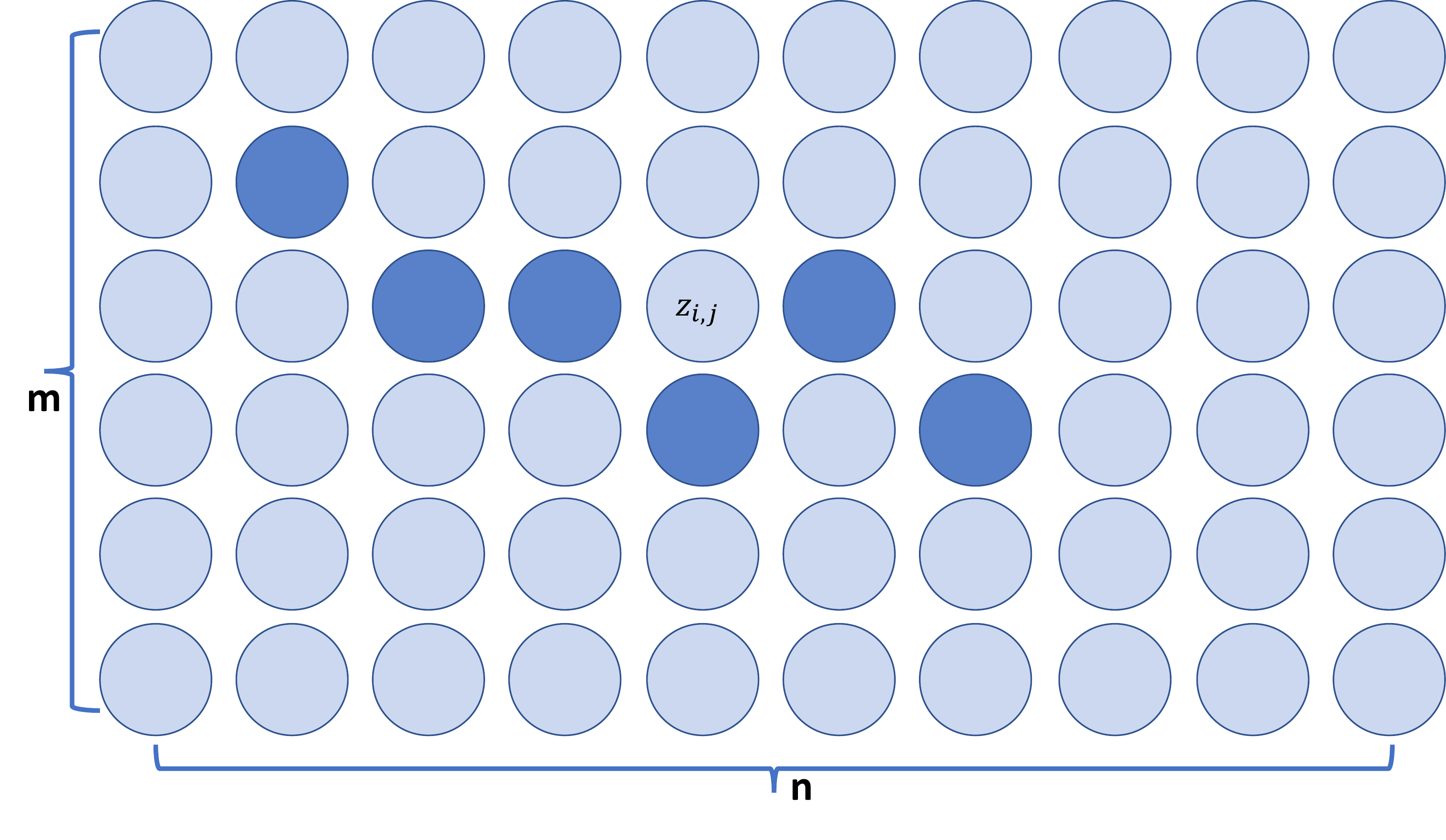}
\caption{The 0-1 index variable $z_{i,j}$ indicates the significance of the corresponding node.}\label{fig:graph2}
    \end{subfigure}
    ~ 
    \begin{subfigure}[b]{0.39\textwidth}
        \includegraphics[width=\textwidth]{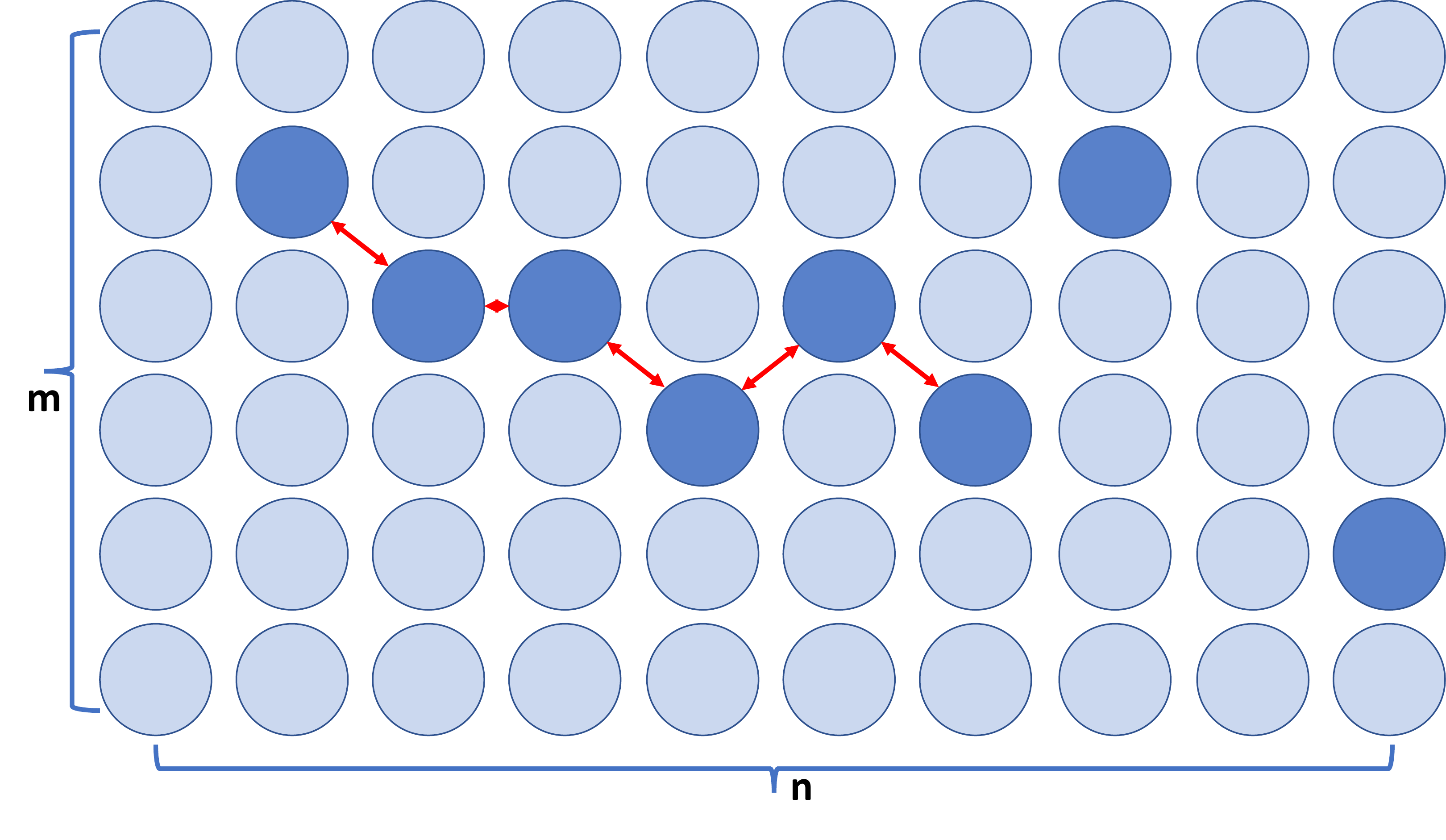}
\caption{The neighborhood relation for significant nodes.}\label{fig:graph3}
    \end{subfigure}
    \caption{Graphical structure for curve detection images}\label{fig:graph}
\end{figure}

In this work, we propose a new strategy, which will not only be able to detect weak signals for the embedded chain, but also have a low order of computational complexity.
In particular, we consider as in the following figure, which explains the graphical structure in an image detection problem.
In Fig. \ref{fig:graph1}, each pixel represents a node in an $m \times n$ graph, with the random variable value $x_{i,j}$, $i = 1, \cdots, m$, $j = 1, \cdots, n$, indicating the intensity of the corresponding pixel.
Upon observing such an image, we first construct an indicator statistic $z_{i,j}$ corresponding to each $x_{i,j}$, which will have the value $z_{i,j} = 1$ if $x_{i,j} > x^\ast$, for some predefined $x^\ast$, and the corresponding node is called a significant node. Otherwise, $z_{i,j} = 0$, and the corresponding node is called an insignificant node.
The significance of nodes is shown in Fig. \ref{fig:graph2}.
We say there is an edge between two significant nodes, if the two nodes are in two consecutive columns and consecutive rows for illustration purpose, which is illustrated in Fig. \ref{fig:graph3}.
Our objective is to detect the existence of the embedded chain.

In the literature, one of the  popular methods for curve detection is using the longest run approach \cite{Aos06Filament, Donoho2005}.
However, when the underlying curve is of length in the order $o(n)$ (instead of $O(n)$), where $n$ is the number of columns in Fig. \ref{fig:graph1}, it will be hard to detect based on only the longest significant run statistic.
While in the literature of change detection, the generalized likelihood ratio statistic (GLR) is a classic way for detecting changes. However, this won't work in our problem due to the unknown locations of the curve.
Furthermore, a complete enumeration of all possible chains with good continuation in an image with $m$-by-$n$ pixels is not a good strategy, because the cardinality of such an enumeration will be approximately $O(e^n)$.

\begin{figure}[htbp]
    \centering
    \begin{subfigure}[t]{0.45\textwidth}
        \includegraphics[width=\textwidth]{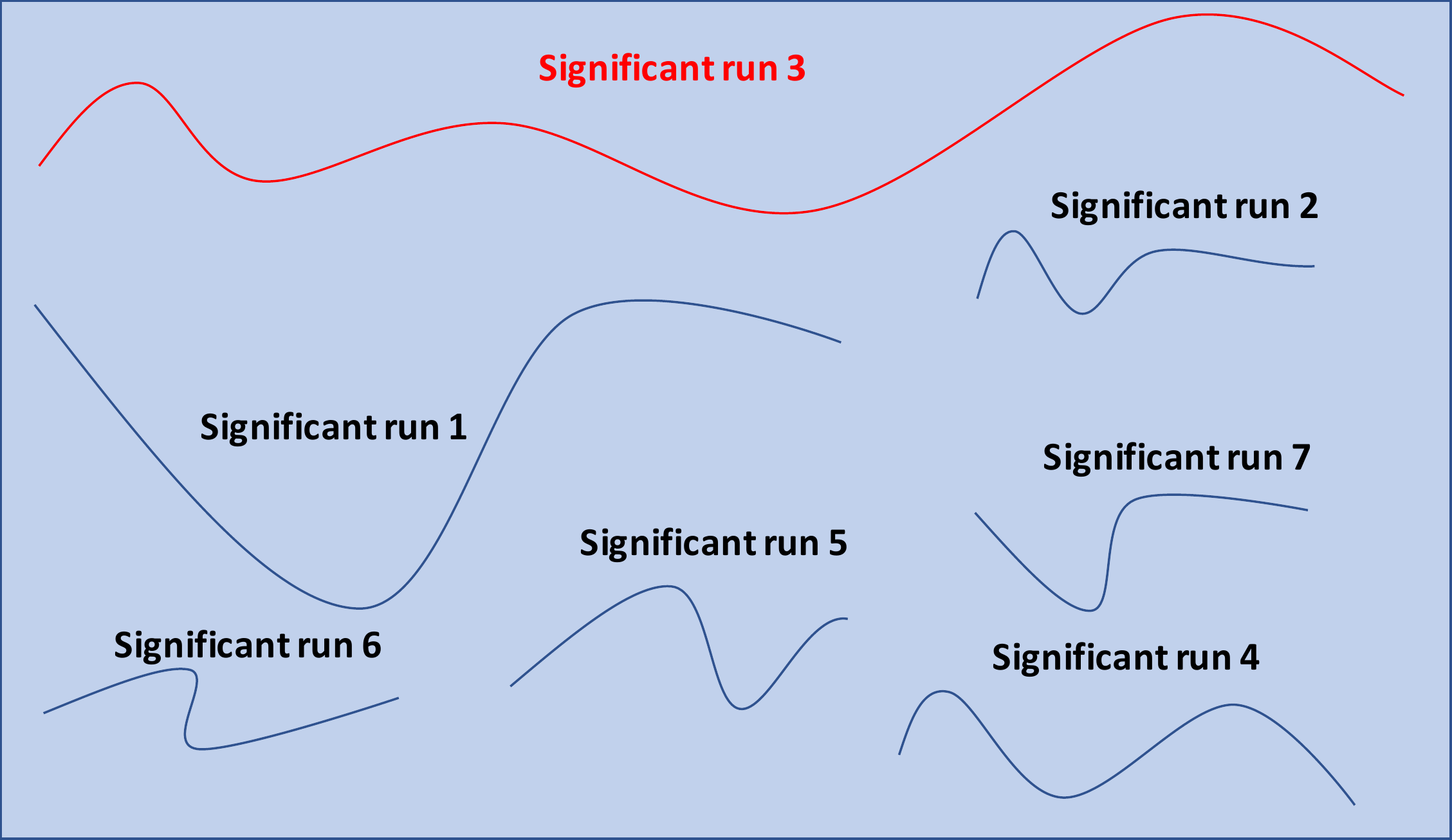}
\caption{Step I: Identify the longest significant run in an image (significant run 3) and accept the existence of an unknown curve if the length exceeds some predefined threshold. Otherwise, go to Step II.}\label{fig:step1}
    \end{subfigure}
    ~ 
    \begin{subfigure}[t]{0.45\textwidth}
        \includegraphics[width=\textwidth]{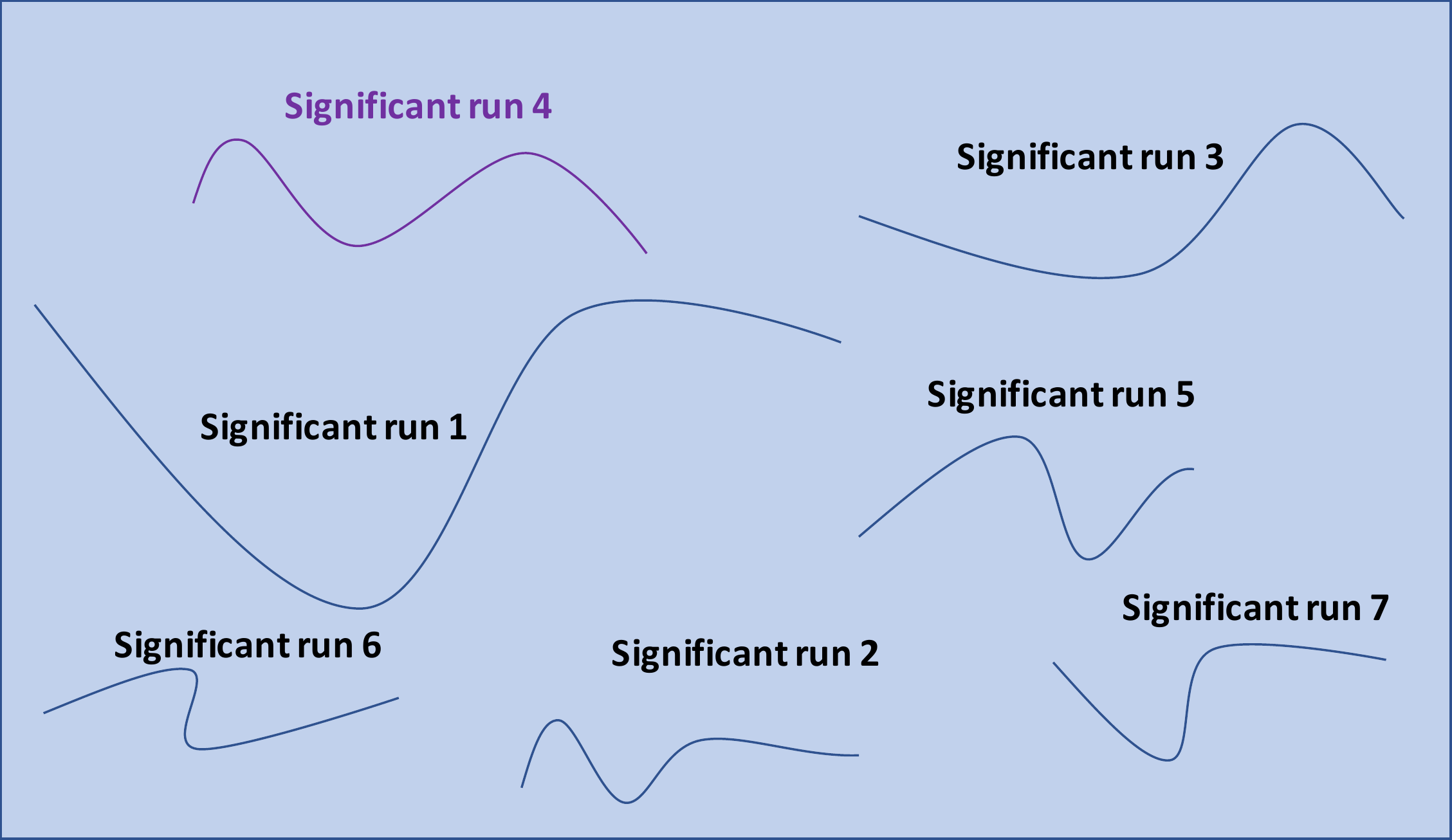}
\caption{Step II: Compute the normalized scan statistics for each significant run and accept the existence of an unknown curve if the maximum scan statistic exceeds some predefined threshold.}\label{fig:step2}
    \end{subfigure}
    \caption{2-step detection algorithm}\label{fig:illustration}
\end{figure}

Instead, we implement a two-step detecting method based on the longest significant chain in an $m$-by-$n$ array as in \cite{JASA06LongSigRun} and the scan statistic for the mean change inspection in the literature of change detection over all the significant runs.
Specifically, we first use a threshold to classify each pixel as either significant or insignificant based on pixel\rq{}s intensity as in Fig. \ref{fig:graph2}, which reduce the problem to a detection problem in a Bernoulli net.
In \cite{JASA06LongSigRun}, the authors derive an asymptotic rate of the length of the longest significant chain in a Bernoulli net.
When applying to our problem with chains consisting of only significant nodes, this technique reduces the number of chains under consideration to a polynomial of $n$.
We identify the longest significant run in the image and  accept the existence of an unknown curve if the length exceeds some predefined threshold as in Fig. \ref{fig:step1}.
Next, in order to detect chains with smaller length, we define a normalized scan statistic to scan the significance over all the significant chains found in the first step as in Fig. \ref{fig:step2}.
We implement our method and find in most cases, the detectable mean lies in between $1$ and $2$.
This result is much better than the detectability by human eyes, in which people can hardly tell the embedded chain with confidence when $\mu<3$ as shown in Fig. \ref{compareFig}.
Our algorithm has three advantages.
First, the algorithm has a very low order of complexity $O(n\log n)$. Note that there are $O(e^n)$ possible chains under consideration, which implies that  our algorithm is very fast.
Second, our detection algorithm is asymptotically powerful, which means as the size of the noisy image becomes larger and larger, the detection errors (type-I error and type-II error) go to zero.
Third, our proposed algorithm is stable. Even if the length of the embedded inhomogeneous chain is as short as $O(\log n)$, the minimum detectable elevated mean in the embedded chain is almost always a half of what can be detected by eyes.
Thus, our algorithm is good in terms of stability.

This paper is organized as follows.
We discuss the statistical model of the detection problem and out proposed algorithm in Section \ref{sec:StatModel}.
Section \ref{sec:LongSigRun}  proves that the type-I error and the type-II error diminish fast, and our algorithm is asymptoticly optimal.
Numeric studies are given in Section  \ref{sec:NumStudy}.
We provide an extension to the case that the number of rows goes to infinity in Section \ref{sec:Extension}.
In order to keep fluency of ideas, the proofs are relegated to the Appendix.

\section{Statistical model}
\label{sec:StatModel}
In this Section, we will first provide the details on our model and the related notations, which will be used throughout this paper. We then illustrate our 2-step algorithm. Related properties of the algorithm are provided in Section \ref{sec:LongSigRun}.

\subsection{Model and Notations}
We consider an $m$-by-$n$ array of nodes $\mathcal{S}$ with $m$ rows and $n$ columns, i.e.,
\begin{equation}\label{defofarray}
\mathcal{S}=\{(i,j): 1\leq i\leq m, 1\leq j\leq n\}.
\end{equation}
An example of $\mathcal{S}$ is shown in Fig. \ref{compareFig} (b) with $m=10$ and $n=60$.

We assume $m\in\mathbb{Z}^{+}$ is fixed and will eliminate this restriction in Section \ref{sec:Extension}. Such an array can be considered as a grid in a two dimensional rectangular region $[1,n]\times[1,m]$.
Assume that each node with coordinate $(i,j)\in\mathcal{S}$ is associated with a normal distributed random variable $X_{i,j}$.
In a digital image, each node indicates a pixel of the image and its corresponding normal random variable $X_{i,j}$ denotes the intensity of the image at $(i,j)$.
We consider all the embedded chains in $\mathcal{S}$ with good continuation such that nodes on the chain are horizontally adjacent and their altitudes are nearly the same---less than $C\in\mathbb{Z}^{+}$ apart for neighboring nodes. To be precise, for $(i,j_0)\in\mathcal{S}$ and $L\geq0$, a chain $\mathcal{L}$ of good continuation with length $L+1$ has the following form,
\begin{eqnarray}\label{defOfChain}
\mathcal{L}&=&\{(i, j_0), (i+1, j_1) \ldots, (i+L, j_L):\nonumber\\
                                                      &&\left|j_k-j_{k-1}\right|\leq C, \text{for } 1\leq k\leq L\}.
\end{eqnarray}

Let $\mathcal{F}_n$ be the set consisting of all chains with good continuation in $\mathcal{S}$.

As a first attempt to formalize matters, consider the problem of testing
\begin{equation}\label{defOfNullHypot}
\mathbb{H}_0: X_{i,j}\stackrel{\mbox{i.i.d.}}{\sim} N(0,\sigma^2), \forall (i,j)\in\mathcal{S}
\end{equation}
versus
\begin{eqnarray}\label{defOfAlterHypot}
\mathbb{H}_1: &X_{i,j}&\stackrel{\mbox{i.i.d.}}{\sim} N(\mu, \sigma^2), \text{ for some }\mu>0,
\mbox{ and an } \\&&\mathcal{L}^0_n \in \mathcal{F}_n, 
\text{ when } (i,j)\in\mathcal{L}^0_n; \label{defofalterhypot}\nonumber \\
&X_{i,j}&\stackrel{\mbox{i.i.d.}}{\sim} N(0,\sigma^2), \forall (i,j)\in\mathcal{S}\setminus\mathcal{L}^0_n,\nonumber
\end{eqnarray}
where $N(\mu, \sigma^2)$ stands for a normal distribution with mean $\mu$ and standard deviation $\sigma$ and $\mathcal{L}^0_n$ is an \textit{unknown} chain in $\mathcal{S}$ with good continuation.
Note that by varying the location and orientation of embedded chain $\mathcal{L}^0_n$ and the value of the parameter $\mu$, there are infinite number of possibilities for the alternative hypothesis $\mathbb{H}_1$ as $n\rightarrow\infty$.
The objective of our forgoing testing problem is to detect whether there \textit{exists} an embedded chain $\mathcal{L}_n^0$ in $\mathcal{S}$ with an elevated mean $\mu>0$ in the digital image $\mathcal{S}$. More specifically, how large the value of $\mu$ and how long the chain $\mathcal{L}^0_n$ should be so that the corresponding alternative hypothesis can be strongly distinguished from the null hypothesis.

Throughout the paper, it is assumed that the length of a chain $\mathcal{L}^0_n$, denoted by $\left|\mathcal{L}^0_n\right|$, can go to $\infty$ as $n\rightarrow\infty$.
We explicitly indicate that $\mathcal{L}^0_n$ depends on $n$ because except Section \ref{sec:Extension}, we assume that $m$ is fixed and our goal it to handle the detectability of (\ref{defOfNullHypot}) versus (\ref{defOfAlterHypot}) asymptotically as the number of columns $n\rightarrow\infty$.
For convenience, in this paper, we assume $\sigma=1$.
Thus our detectability problem is equivalent to finding the SNR (signal-to-noise ratio $|\frac{\mu}{\sigma}|$) threshold, which ensures the asymptotic powerfulness of the hypothesis testing problem.
The length of the chain and the number of nodes on the chain are the same.
We use $\left|\cdot\right|$ to denote the length of a chain or the cardinality of a set. We use $C, C_1, C_2, \delta_1, \delta_2, \eta, \eta_1, \eta_2, \zeta$ to indicate positive constants which may vary case by case.

\subsection{Algorithm}
\label{sec:Summary}
Given the notations in the previous section, the following shows our detection problem when $\sigma = 1$:
\[
\mathbb{H}_0: X(i,j)\sim N(0,1), i.i.d., \forall(i,j)\in\mathcal{S}
\]
versus

\begin{eqnarray}
\nonumber
\mathbb{H}_1: &X(i,j)&\sim\mu+N(0,1), i.i.d., \forall(i,j)\in\mathcal{L}^{0}_{n}, \mbox{ for an }\\ && \mathcal{L}^0_n \in \mathcal{F}_n, \mbox{ and some } \mu >0,
\nonumber \\
&X_{i,j}&\stackrel{\mbox{i.i.d.}}{\sim} N(0,1), \forall (i,j)\in\mathcal{S}\setminus\mathcal{L}^0_n,\nonumber
\end{eqnarray}

Given a prescribed threshold $x^{\ast}$, we say a node $(i,j)\in\mathcal{S}$ is significant, if its corresponding observed pixel value satisfies $X(i,j)>x^{\ast}$. Recall that
\[
\mathcal{S}=\{(i,j): 1\leq i\leq m, 1\leq j\leq n\}
\]
is the set of pixels in the image.
We define a random variable $z:\mathcal{S}\rightarrow\{0,1\}$ to indicate the significance of nodes, i.e., $z(i,j)=1$ if $X(i,j)>x^{\ast}$ and $z(i,j)=0$ otherwise.
It is easy to see that under the null hypothesis $\mathbb{H}_0$,  $z(i,j)$ are i.i.d. Bernoulli distributed with success probability
\[
p=\mathbb{P}(N(0,1)>x^{\ast}),
\]
where we use the expression $\mathbb{P}(N(0,1)>x^{\ast})$ to denote the probability of exceeding $x^{\ast}$ of a random variable $X\sim N(0,1)$. Similar expressions are used for different random variables in the context without misunderstanding.
Given this notation, a chain with good continuation $\mathcal{L}\in\mathcal{F}_n$ is said to be significant, if all the nodes on the chain are significant.
Define a random variable
\begin{equation}\label{defOfZ}
Z:\mathcal{F}_n\rightarrow\{0,1\}
\end{equation}
such that $Z(\mathcal{L})=1$ if $\mathcal{L}\in\mathcal{F}_n$ is significant and $Z(\mathcal{L})=0$ otherwise.
Obviously, the probability of a chain $\mathcal{L}$ to be significant under $\mathbb{H}_0$ is the following,
\begin{equation}\label{ProbOfSigRun}
\mathbb{P}(Z(\mathcal{L})=1)=p^{\left|\mathcal{L}\right|}.
\end{equation}

Given the aforementioned notations (under $\mathbb{H}_0$), let us recall a series of results in \cite{JASA06LongSigRun}, which states the asymptotic rate of the length of the longest significant chain with good continuation in the $m$-by-$n$ array of nodes as $n\rightarrow\infty$.
Throughout the paper, let $L_0(n)$ be the longest significant chain in the image and let $\left|L_0(n)\right|$ be its length.
Notice that $\left|L_0(n)\right|$ actually depends on $m, C$ and $p$ in addition to $n$, but we simplify the notation because all the parameters except $n$ are constants.
Let the following
\[
P_n=\mathbb{P}_{m,C,p}(\left|L_0(n)\right|=n)
\]
denote the probability that the length of the longest significant run is $n$, when there are exactly $n$ columns.
The following lemma in \cite{JASA06LongSigRun} states the asymptotic ratio of $P_n/P_{n-1}$.
\begin{lemma}\label{lemmaRatio}
Define $\rho_n=P_n/P_{n-1}$.
There exists a constant $\rho\in (0,1)$ that depends only on $m, C, $ and $p$, but not on $n$, such that
\begin{equation}\label{RatioOfAcross}
\lim_{n\to\infty}\rho_n=\rho.
\end{equation}
\end{lemma}

\begin{remark}
We say a significant chain is across if and only if it traverses all columns.
The ratio $\rho_n$ is the conditional probability that there is an across significant chain for $n$ columns, conditioning on the fact that there is an across significant run in the previous $(n-1)$ columns.
We may call this the chance of preserving across significant chains.
Lemma \ref{lemmaRatio} shows that as the number of columns goes to infinity, the chance of preserving across significant chains converges to a constant. Table \ref{tbRho} gives the exact values of $\rho$ for different $p$'s and $m$'s: $m=4, 8, 10$ and $C=1$.
As one can expect, $\rho$ increases as $p$ increases. The algorithmic complexity for finding $\rho$ can be $O(2^{3m})$ \cite{JASA06LongSigRun}.
\end{remark}
\begin{table}[!htbp]
\caption{The values of $\rho$ for different values of $m$ and $p$,
when $C=1$.}
\begin{center}
\begin{tabular}{|c|c|c|c|c|c|c|}
\hline
$\rho$&p = 0.1&p = 0.2&p = 0.3&p = 0.4&p = 0.5 & p = 0.6\\
\hline
m=4&0.2444 &   0.4564& 0.6341 & 0.7758& 0.8804 & 0.9482\\
m=8& 0.2654 & 0.4955& 0.6869 & 0.8363 & 0.9383& 0.9876\\
m=10&0.2691 & 0.5022& 0.6958& 0.8467& 0.9486& 0.9930\\
\hline
\end{tabular}
\end{center}
\label{tbRho}
\end{table}

We are now ready to describe our algorithms and the detectability conditions, which can separate the null hypothesis $\mathbb{H}_0$ from the alternative ones $\mathbb{H}_1$, in the following.
 This test is independently considered in a series of papers such as \cite{perim4} and \cite{perim5}.

One of our focus is on the test that rejects for large values of the following scan statistic:
\begin{equation}\label{eq:scan}
\max_{\mathcal{L}\in\mathcal{F}_n}\sum_{(i,j)\in\mathcal{L}}\frac{X(i,j)}{\sqrt{\left|\mathcal{L}\right|}}.
\end{equation}
The normalization in (\ref{eq:scan}) is such that each term in the maximization is standard normal under the null hypothesis and allows us to compare chains of different sizes.
The scan statistic was originally proposed in the context of cluster detection in point clouds in \cite{scanStat}.
The scan statistic is the prevalent method in disease outbreak detection, with many variations (See \cite{Evaluation} and \cite{treeBasedScan}).
In our paper, $\mathcal{F}_n$ is the set consisting of all the chains with good continuation in the image.

However, it is not hard to see that there are $O(e^n)$ chains in $\mathcal{F}_n$.
We thus will not use the scan statistics directly, but rather restrict the scanning to a subset of $\mathcal{F}_n$ as follows.

Recall that $Z: \mathcal{F}_n\rightarrow\{0,1\}$, which is defined in (\ref{defOfZ}), is an indicator of significance for all chains with good continuation (i.e., in $\mathcal{F}_n$). Let $\mathcal{E}_n$ be a random subset of $\mathcal{F}_n$ consisting of all significant chains, i.e.,
\begin{equation}\label{defOfRanSigSet}
\mathcal{E}_n=\{\mathcal{L}\in\mathcal{F}_n: Z(\mathcal{L})=1\}.
\end{equation}

Recall that the longest significant chain in $\mathcal{E}_n$ is denoted by $L_0(n)$ and that $\left|L_0(n)\right|$ is the length of $L_0(n)$.
The next theorem in \cite{JASA06LongSigRun} provides the asymptotic rate of $\left|L_0(n)\right|$, which is a generalization of the well-known Erd$\ddot{o}$s-R$\acute{e}$nyi law (See \cite{LofLHR, PetrovVV, NLawofLN}).
\begin{theorem}\label{thConvRate}
Under the null hypothesis, as $n\rightarrow\infty$,
\begin{equation}\label{eqConvRate}
\frac{\left|L_0(n)\right|}{\log_{1/\rho}n}\rightarrow 1, \text{\quad almost surely}.
\end{equation}
A special case of the theorem is that when $m=1$, then (\ref{eqConvRate}) holds with $\rho$ replaced by $p$.
\end{theorem}

Given this result, it is easy to obtain the following observation, which states the relation of $\rho$ and $(m, p)$.
Since $\left|L_{0}(n)\right|$ actually depends on $p$ and $m$, we use the notation $\left|L_{0}(m,n,p)\right|$ in the next corollary to make this dependence explicit.

\begin{corollary}\label{co:depofrho}
Given a pair of positive integers $m_1,m_2$ and a pair of probabilities $p_1, p_2$ with $m_1\leq m_2$ and $p_1\leq p_2$, we have
\[
\rho(m_1,p_1)\leq\rho(m_2,p_1)\text{\quad and \quad}\rho(m_1,p_1)\leq \rho(m_1,p_2)
\]
\end{corollary}

 By Theorem \ref{thConvRate} and Egoroff's Theorem (See \cite{Royden2010}), given any small $\epsilon>0$ and $\delta>0$, there exists a large $N\in\mathbb{Z}^{+}$, such that for all $n\geq N$ with probability $1-\delta$ under $\mathbb{H}_0$, we have
\begin{equation}\label{eqEgoroff}
\left|\frac{\left|L_0(n)\right|}{\log_{1/\rho}n}-1\right|<\epsilon.
\end{equation}
Let $b_{\epsilon,n} = (1+\epsilon)\log_{1/\rho}n$, which under $\mathbb{H}_0$, is the upper bound of the length of the longest significant run with probability $1-\delta$.

It is not hard to see that under the null hypothesis $\mathbb{H}_0$, for large $n\geq N$, with probability $1-\delta$, an upper bound of $\left|\mathcal{E}_n\right|$ under $H_0$ is
\begin{equation}\label{upperBoundOfSigRun}
\sum_{k=1}^{b_{\epsilon,n}+1}mn(2C+1)^{k-1}=mn[(2C+1)^{b_{\epsilon,n}+1}-1]/2C.
\end{equation}

By Corollary \ref{co:depofrho}, we can choose the threshold of significance, $x^{\ast}$, large enough, so that
\[
p=\mathbb{P}(N(0,1)>x^{\ast})
\]
is small enough to render $\rho<\frac{1}{2C+1}$. By Corollary \ref{co:depofrho}, it is easy to see that $p\leq\rho<\frac{1}{2C+1}$. Thus an upper bound on $\left|\mathcal{E}_n\right|$, according to Equation (\ref{upperBoundOfSigRun}) and the definition of $b_{\epsilon,n}$, is
\begin{eqnarray}
&&mn\frac{2C+1}{2C}\cdot(2C+1)^{\log_{2C+1}n\cdot(1+\epsilon)\log_{1/\rho}(2C+1)}\nonumber\\
&\leq& mn\frac{2C+1}{2C}n^{(1+\epsilon)\log_{1/\rho}(2C+1)}.\label{FurtherUpperBound}
\end{eqnarray}
Since $\epsilon$ is arbitrary, as $n$ becomes large, the right hand side of (\ref{FurtherUpperBound}) is less than $mn^{2-\delta_0}$ for some $\delta_0>0$.
Thus under the null hypothesis $\mathbb{H}_0$, $\left|\mathcal{E}_n\right|$ grows slower than $mn^2$.

Now, we define our statistic based on all significant chains in the array $\mathcal{S}$.

\begin{definition}\label{defOfAllSigChain}
In an array of $m$-by-$n$ nodes $\mathcal{S}$, let $X(i,j)$ be the normally distributed random variable associated with each node $(i,j)$. Let $X(\mathcal{L}) = \sum_{(i,j)\in\mathcal{L}}\frac{X(i,j)}{\sqrt{\left|\mathcal{L}\right|}}.$
Then we define a statistic based on all significant chains to be
\begin{equation}\label{defOfCompSigStat}
X_{s}^{\ast}=\max_{\mathcal{L}\in\mathcal{E}_n}X(\mathcal{L}).
\end{equation}
\end{definition}

Given the aforementioned notations of $x^{\ast}$, $\mathcal{E}_n$ and $\rho$, we now describe the algorithm for the analysis of a noisy image $\mathcal{S}=\{X(i,j), 1\leq i\leq m, 1\leq j\leq n\}$, looking for the existence of suspected chains with good continuation in $\mathcal{S}$.
The algorithm has two steps and its computational complexity is $O(n\log n)$.
\begin{itemize}
\item \textbf{Step I:} Count the length of the longest chain $L_0(n)$ in $\mathcal{E}_n$. If the length
\[
\left|L_0(n)\right|>(1+\epsilon/2)\log_{1/\rho}n
\]
for some small $\epsilon>0$, then reject $\mathbb{H}_0$; otherwise, go to Step II.
\item \textbf{Step II:} Compute $X_s^{\ast}$ as in (\ref{defOfCompSigStat}). If
\[
X^{\ast}_s>\sqrt{2(1+\delta_2)\log n},
\]
for some small $\delta_2>0$, then reject $\mathbb{H}_0$; otherwise, accept $\mathbb{H}_0$.
\end{itemize}

\subsection{Computation cost}

In this section, we will show that the first step takes $O(n)$ flops and the second $O(n\log n)$ using the dynamic programming approach. Hence this algorithm takes $O(n\log n)$ flops in total with $O(n\log n)$ required space for storage.

\textbf{Algorithm to find $\left|L_0(n)\right|$:} Recall that for a node $(i,j)\in\mathcal{S}$, we use $z(i,j)=1$ $(=0)$ to denote the significance (insignificance) of $(i,j)$. Given a realization $\{X(i,j): 1\leq i\leq m, 1\leq j\leq n\}$, let $Y_1$ be $\{Y_1(i,j): 1\leq i\leq m, 1\leq j\leq n\}$ such that
\begin{eqnarray*}
Y_1(i,1)&=&z(i,1), \text{for } i=1,\ldots, m;\\
Y_1(i,j)&=&z(i,j)[1+\max_{i'\in\Omega(i)}Y_1(i', j-1)], \\
          &\text{for }& i=1,\ldots, m, j=2,\ldots,n,
\end{eqnarray*}
where $\Omega(i)=\{i': \left|i'-i\right|\leq C, 1\leq i'\leq m\}$ denotes the set containing neighboring indices of $i$. Finally, the value  $\left|L_0(n)\right|$ can be computed as follows:
\[
\max_{(i,j)\in\mathcal{S}}Y_1(i,j).
\]
It is not hard to see that this algorithm takes $Cmn$ time for $C>0$.

\textbf{Algorithm to find $X_s^{\ast}$:} We use $z(i,j)=1 (=0)$ to denote the (in)significance of $(i,j)\in\mathcal{S}$. Let $U=\lceil3\log_{1/\rho}n\rceil$ and $Y_2=\{Y_2(i,j,u):1\leq i\leq m, 1\leq j\leq n, 1\leq u\leq U\}$ such that
\[   Y_2(i,j,u) = \begin{cases}
      X(i,j), & u = 1;\\
            Y_2(i,j,u-1), &  i=1,\ldots, m, \\
            &j=1,\ldots, u-1, \\
            &2\leq u\leq U;\\
      z(i,j)[X(i,j)+&  i=1\ldots, m,\\
      \max_{i'\in\Omega(i)}Y_2(i',j-1,& j=u\ldots, n,\\
      u-1)],& 2\leq u\leq U; \\
\end{cases}\]

 and
 \begin{eqnarray*}
X^{\ast}_s&=&\max_{(i,j)\in\mathcal{S}, 1\leq u\leq U}\frac{Y_2(i,j,u)}{\sqrt{u}},
\end{eqnarray*}
where $\Omega(i)$ is the set of neighboring indices of $i$ as in the algorithm to find $\left|L_0(n)\right|$.
It is easy to see that this algorithm takes $Cmn\log n$ time for some $C>0$.
\begin{remark}
Under the alternative hypothesis $\mathbb{H}_1$, if the information on the length of the unknown chain $\mathcal{L}^0_n$ is available, which is in the order of $O(n^\alpha)$ for some $0<\alpha<1$, by conducting only Step I in the previous algorithm, we will be able to make an optimal decision.
The corresponding testing problem is asymptotically powerful given the signal $\mu \geq \mu^\ast$, where $\mu^\ast$ is such that $p_1 > \rho^{\frac{\alpha}{1+\epsilon}}$ for some small $\epsilon > 0$.
In this scenario, the computation complexity is just $O(n)$.
\end{remark}

\section{Asymptotic properties}
\label{sec:LongSigRun}

In this section, we give the upper bounds for type-I error under $\mathbb{H}_0$ and type-II error under $\mathbb{H}_1$. We also prove the asymptotic optimality of our proposed detection algorithm.

\subsection{Behavior under $\mathbb{H}_0$}
We need to show that with overwhelming probability, under $\mathbb{H}_0$, there will be no significant chain longer than $(1+\epsilon)\log_{1/\rho}n$ and $X_s^{\ast}<\sqrt{2(1+\delta_2)\log n}$, for any small $\epsilon>0$ and $\delta_2>0$, as $n\rightarrow\infty$.
The former is shown in (\ref{eqEgoroff}) and we show the latter in the following theorem.

\begin{theorem}\label{type1error}
Under the null hypothesis $\mathbb{H}_0$, for any small $\delta>0$, there exists a constant $\sigma_1$, depending on $p$ and a large $N\in\mathbb{Z}^{+}$, such that for any $n\geq N$, we have the following:
\begin{equation}\label{upperBoundfirstErr}
\mathbb{P}(X_s^{\ast}>\tau) \leq mn\sigma_1\exp\left\{-\frac{\tau^2}{2}\right\}+\delta.
\end{equation}
Thus for any $\delta_2>0$, when $\tau = \tau_s^{\ast}=\sqrt{2(1+\delta_2)\log mn}$, we have $\mathbb{P}_{\mathbb{H}_0}(X_{s}^{\ast}>\tau_s^{\ast})\rightarrow0$ as $n\rightarrow\infty$.
\end{theorem}

So far, we have studied the asymptotic behavior of $X^{\ast}_s$ under $\mathbb{H}_0$ and prove that the type-I error tends to $0$ as $n\rightarrow\infty$.
In the next subsection, we delve into the behavior of $X_s^{\ast}$ under $\mathbb{H}_1$ and shall prove the diminishing type-II error.

\subsection{Asymptotic behavior under $\mathbb{H}_1$}\label{sebsec:asymH1}
We first show the condition under which the type-II error diminishes fast under $\mathbb{H}_1$, if the underlying chain with elevated mean satisfies that $\left|\mathcal{L}^0_n\right|\geq\zeta_1 n$ for some $\zeta_1>0$. Let us first see the behavior of the longest significant chain embedded in $\mathcal{L}^0_n$ under a specific alternative hypothesis: $\mathbb{H}_1(\mathcal{L}^0_n,\mu)$.  Denote $\mathbb{P}(N(\mu,1)>x^{\ast})$ by $p_1$, which is the probability of nodes to be significant in the chain $\mathcal{L}^0_n$.
Let $L_1(n)$ be the longest significant chain in $\mathcal{L}^0_n$ and let $\left|L_1(n)\right|$ be its length.
Recall that $\left|L_0(n)\right|$ is the length of the longest significant chain in the image and so $\left|L_1(n)\right|\leq\left|L_0(n)\right|$ since $L_1(n)\in\mathcal{E}_n$.
By a special case of Theorem \ref{thConvRate} ($\rho=p_1$ when $m=1$), It is obvious that as
\[
\left|\mathcal{L}^0_n\right|\rightarrow\infty,
\]
we have the following convergence rate of $\left|L_1(n)\right|$,
\begin{equation}\label{extractLongSigRun}
\frac{\left|L_1(n)\right|}{\log_{1/p_1}\left|\mathcal{L}^0_n\right|}\rightarrow 1\text{\quad almost surely}.
\end{equation}

Therefore by Ergoroff's theorem, given any small $\epsilon_1>0$ and $\delta>0$, there exists a large $N\in\mathbb{Z}^{+}$ such that for all $n\geq N$ with probability $1-\delta$ we have
\begin{equation}\label{egoroffunderH1}
\left|\frac{\left|L_1(n)\right|}{\log_{1/p_1}\left|\mathcal{L}^0_n\right|}-1\right|<\epsilon.
\end{equation}
That is to say,  with probability at least $1-\delta$, one has
\[
(1-\epsilon_1)\log_{1/p_1}\left|\mathcal{L}^0_n\right|\leq \left|L_1(n)\right|\leq(1+\epsilon_1)\log_{1/p_1}\left|\mathcal{L}^0_n\right|.
\]
We will prove the following theorem regarding the type-II error.

\begin{theorem}\label{thLongSigChain}
In an array $\mathcal{S}$ of $m$-by-$n$ nodes, consider the following detection problem
\[
\mathbb{H}_0: X(i,j)\sim N(0,1), i.i.d., \forall(i,j)\in\mathcal{S}
\]
versus
\begin{eqnarray*}
\mathbb{H}_1: &&X(i,j)\sim\mu+N(0,1), i.i.d., \forall(i,j)\in\mathcal{L}^{0}_{n}, \\
&&\mbox{ for some } \mathcal{L}^0_n \in \mathcal{F}_n \mbox{ with  } \left|\mathcal{L}^0_n\right| \geq \zeta_1 n^\alpha, \mbox{ and some   } \mu >0,
\end{eqnarray*}
where $\mathcal{L}^0_n$ is a chain with good continuation ($C$ apart) with length in the order of $O(\zeta_1 n^\alpha)$ for some $\zeta_1>0$ and $0<\alpha \leq1$. If $\mu$ is such that
\begin{equation}\label{condOfLongSigRun}
p_1>\rho^{\alpha\log \zeta_1 n/(1+\epsilon)\log n}\rightarrow \rho^{\alpha/(1+\epsilon)}\text{ as } n\rightarrow\infty
\end{equation}
for some small $\epsilon>0$, then as $n\rightarrow\infty$, we have
\begin{equation}\label{asympPowerOfLongSigRun}
\mathbb{P}(\left|L_0(n)\right|>(1+\epsilon/2)\log_{1/\rho}n\big|\mathbb{H}_1)\rightarrow 1.
\end{equation}
\end{theorem}
Let us consider the case under $\mathbb{H}_1$, where the assumptions on $\left|\mathcal{L}^0_n\right|$ in Theorem \ref{thLongSigChain} fail.
Denote two constants $\lceil(1-\epsilon)\log_{1/p_1}\left|\mathcal{L}^0_n\right|\rceil$ and $\lfloor(1+\epsilon)\log_{1/p_1}\left|\mathcal{L}^0_n\right|\rfloor$ by $c^{L}_{\epsilon,n}$ and $c^{U}_{\epsilon,n}$ respectively.
Recall a chain $\mathcal{L}$ is in $\mathcal{E}_n$ if and only if $X(i,j)>x^{\ast}$ for every node $(i,j)\in\mathcal{L}$, where $x^{\ast}$ is the threshold of significance.
In Definition \ref{defOfAllSigChain}, $X_s^{\ast}$ is said to be the maximum of all $X(\mathcal{L})$ among $\mathcal{L}\in\mathcal{E}_n$, which is of course no smaller than $X(L_1(n))$.
We now give the asymptotic diminishing rate of the type-II error.

\begin{theorem}\label{type2error}
Under the alternative hypothesis $\mathbb{H}_1$, for any small $\delta>0$ and $\epsilon>0$, if $\tau_s^{\ast} = \sqrt{2(1+\delta_2)\log mn} < \mu\sqrt{c^{L}_{\epsilon,n}}$, which is defined in Theorem \ref{type1error}, we have the following:
\begin{equation*}
\mathbb{P}(X^{\ast}_s>\tau_s^{\ast} \big|\mathbb{H}_1)\rightarrow1, \text{\quad as\quad} n\rightarrow\infty.
\end{equation*}
\end{theorem}

Since $\epsilon$ in the aforementioned theorem is arbitrary, $c^{L}_{\epsilon,n}\approx \log_{1/p_1}\left|\mathcal{L}^0_n\right|$, we may change the condition $\tau_s^{\ast} < \mu\sqrt{c^{L}_{\epsilon,n}}$ in Theorem \ref{type2error} to $\tau_s^{\ast} < \mu\sqrt{\log_{1/p_1}\left|\mathcal{L}^0_n\right|}$.
Theorem \ref{type2error} is also applicable in the case $\left|\mathcal{L}^0_n\right| = O(n^\alpha)$.

\subsection{Asymptotical optimality}
\begin{theorem}
Under the assumptions of Theorem \ref{thLongSigChain}, \ref{type2error}, as $n\rightarrow\infty$, we have
\begin{equation}\label{simConvTo0}
\mathbb{P}(\text{accept }\mathbb{H}_1\big|\mathbb{H}_0)+\mathbb{P}(\text{accept }\mathbb{H}_0\big|\mathbb{H}_1)\rightarrow 0,
\end{equation}
which means our algorithm is asymptotically powerful in terms of the following definition in \cite{Donoho2005}.
\end{theorem}

\begin{definition}\label{asymptoticPower}
In a sequence of testing problems $(\mathbb{H}_{0,n})$ versus $(\mathbb{H}_{1,n})$, we say that a sequence of tests $T_n$ is asymptotically powerful if
\[
\mathbb{P}_{\mathbb{H}_{0,n}}\{T_n\text{ rejects } \mathbb{H}_{0,n}\}+\mathbb{P}_{\mathbb{H}_{1,n}}\{T_n\text{ accepts } \mathbb{H}_{0,n}\}\rightarrow 0,
\]
as $n\rightarrow\infty$.
\end{definition}

\section{Numerical study}
\label{sec:NumStudy}
In this section, we carry out  numerical studies on the detection problem of the suspected chain with good continuation in a noisy image.
We will explain our detectability using simulations, which show the minimal mean required for the chain to be detectable.
For simplicity, in all the simulated examples in Subsections \ref{subsec:On} and \ref{subsec:smallon}, it is assumed to have $C=1$ and $m=10$.
We then use the well-known solar flare example to show how our proposed procedures can be used for the detection of solar flare in Subsection \ref{subsec:solar}.

\subsection{$\left|\mathcal{L}_n^0\right|\sim O(n^\alpha)$}\label{subsec:On}
\begin{figure}[hbt]
\centering\includegraphics[width=90mm,height = 10mm]{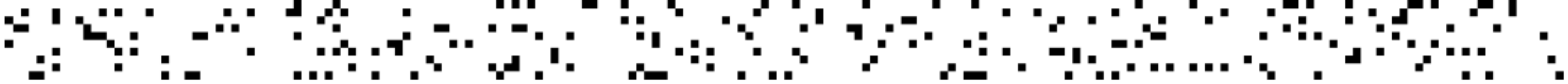}
\caption{Black nodes are significant while white nodes are not significant.}\label{sigNodes}.
\end{figure}

In this subsection, we assume that $\left|\mathcal{L}_n^{0}\right|\geq \zeta_1\cdot n$ for some unknown constant $1>\zeta_1>0$.
From Table \ref{tbRho}, we can choose $x^{\ast}$ to be the $90th$ percentile of the standard normal distribution. Thus $x^{\ast}=1.2816$ and $\rho=0.2691<\frac{1}{3}$.
In Fig. \ref{sigNodes} all the significant nodes are black i.e., $\{(i,j)\in\mathcal{S}: X(i,j)>x^{\ast}\}$, while non-significant nodes are white.

Let $\epsilon=0.0001$ and thus as shown in Inequality (\ref{condOfLongSigRun}), $\mu$ should satisfy
\begin{eqnarray}
p_1&=&\mathbb{P}(N(\mu,1)>x^{\ast})\nonumber\\
&=&\mathbb{P}(N(\mu,1)>1.2816)\nonumber\\
&>&\rho^{\log\zeta_1 n/1.0001\log n}\label{conditionOfMu}
\end{eqnarray}
Let $x_q$ be the $q$th percentile for the standard normal distribution, i.e., $q/100 = \mathbb{P}(N(0,1)>x_q)$.
Let $p_1^\rho = \rho^{\log\zeta_1 n/1.0001\log n}$, $q_1^\rho = 1 - p_1^\rho$.
Thus (\ref{conditionOfMu}) is equivalent to find: $\mu > 1.2816 - x_{100q_1^\rho}$.
Table \ref{minimumOfMu} gives value of $\mu$ which satisfies (\ref{conditionOfMu}) according to different values of $\zeta_1$ and $n$.
As we can see when $\zeta_1$ or the number of columns $n$ increase, we have smaller value of $\mu$ in the table, which indicates stronger detectability in the noisy image.
Intuitively, the increasing length of the inhomogeneous chain under $\mathbb{H}_1(\mathcal{L}^0_n,\mu)$ yields strong visibility.

\begin{table}[!htbp]
\caption{The minimum detectability of $\mu$ when $C=1$ and $m=10$ and $\left|\mathcal{L}_n^0\right|=\zeta_1 n$.}
\begin{center}
\begin{tabular}{|c|c|c|c|c|c|c|}
\hline
$\zeta_1$ & 1/10 & 1/5 & 1/4 & 1/3 & 1/2 & 1\\
\hline
n=$2\times10^2$ & 1.2216 & 1.0307 & 0.9745 & 0.9052 & 0.8126 & 0.6661\\
n=$3\times10^2$ & 1.1740 & 1.0017 & 0.9504 & 0.8869 & 0.8017 & 0.6661\\
n=$5\times10^2$ & 1.1247 & 0.9710 & 0.9249 & 0.8675 & 0.7901 & 0.6661\\
n=$10^3$& 1.0716 & 0.9375 & 0.8969 & 0.8461 & 0.7772 & 0.6661\\
n=$2\times10^3$ & 1.0296 & 0.9107 & 0.8743 & 0.8288 & 0.7668 & 0.6661\\
n=$5\times10^3$ & 0.9860 & 0.8824 & 0.8506 & 0.8105 & 0.7556 & 0.6661\\
n=$10^4$ & 0.9594 & 0.8650 & 0.8359 & 0.7991 & 0.7487 & 0.6661\\
n=$10^5$ & 0.8960 & 0.8232 & 0.8004 & 0.7716 & 0.7319 & 0.6661\\
n=$10^6$ & 0.8553 & 0.7959 & 0.7772 & 0.7535 & 0.7207 & 0.6661\\
\hline
\end{tabular}
\end{center}
\label{minimumOfMu}
\end{table}

\begin{figure}[!htbp]
\centering\includegraphics[width=100mm]{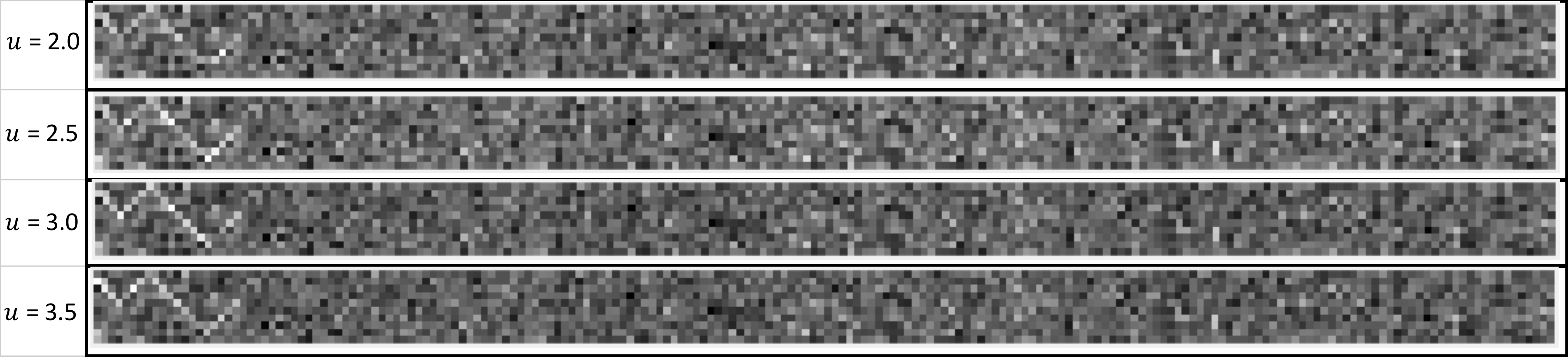}
\caption{Grayscale images of $10\times200$ pixels with different means under $\mathbb{H}_1$ for a chain of length $20$. When the elevated mean is less than $2.5$, it is very hard to identify the inhomogeneous chain.}\label{Fi:n200c01}.
\end{figure}

\begin{figure}[!htbp]
\centering\includegraphics[width=100mm]{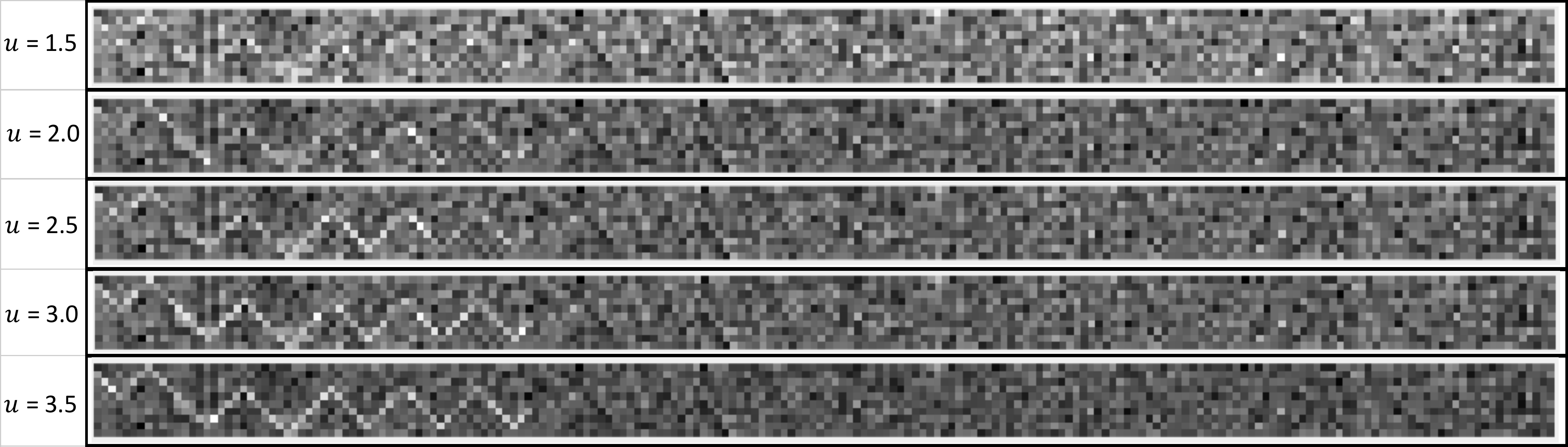}
\caption{Gray-scale images of $10\times300$ pixels with different means under $\mathbb{H}_1$ for a chain of length $60$. The inhomogeneous chain with good continuation becomes apparent when $\mu=2.5$.}\label{Fi:n300c02}.
\end{figure}

Below is the simulation result for $m=10$, $n=200$, $C=1$ and $\zeta_1=\frac{1}{10}$. In Fig. \ref{Fi:n200c01}, when $\mu\leq 2.5$, by human eyes, it is hard to tell whether there is an embedded chain different from the background. However, our method works for $\mu>1.2216$.
Fig. \ref{Fi:n300c02} gives a simulation for $m=10$, $n=300$, $C=1$ and with $\frac{1}{5}$ portion of nodes on a chain with good continuation.
We can see in Fig. \ref{Fi:n300c02} the chain becomes apparent when $\mu\geq2.5$ and our theory supports the detectability of such a chain when $\mu>1.1740$.

\subsection{$\zeta\log n<\left|\mathcal{L}^0_n\right|<c n^{1-\delta}$}\label{subsec:smallon}
In this part, we consider the minimum detectability when $\left|\mathcal{L}_n^0\right|$ is of order $o(n)$, i.e.,
\[
\zeta\log n<\left|\mathcal{L}^0_n\right|\leq c\cdot n^{1-\delta}
\]
for some $\delta>0$ such as $\left|\mathcal{L}_n^0\right|=c_1\sqrt{n}$ and $\left|\mathcal{L}^0_n\right|=c_2\log n$, where $c_2>\zeta$. In both cases, the value of the elevated mean $\mu$ that can be detectable is within our expected range.

\subsubsection{$\left|\mathcal{L}^0_n\right|=c\sqrt{n}$}
Similar as in Section \ref{subsec:On}, we also provide the value of $\mu$ which satisfies  the inequality  (\ref{condOfLongSigRun}), according to different values of $\zeta_1$ and $n$.
The minimum value of $\mu$ that satisfies (\ref{condOfLongSigRun}) is listed in Table \ref{OnSqrtTable}.
Moreover, based on the results in Theorem \ref{type2error}, we can also compute another set of detectability threshold for the mean signal, where the minimum value of $\mu$ should  satisfy (\ref{exMu}):
\begin{equation}\label{exMu}
\mu\sqrt{\log_{1/p_1}(c\sqrt{n})}>\sqrt{(2+\delta_2)\log mn}.
\end{equation}
However, this set of detectability is not as tight as using Theorem \ref{thLongSigChain}.

\begin{table}[!htbp]
\caption{The minimum detectability of $\mu$ when $m=10$, $C=1$ and $|\mathcal{L}^0_n|=c\sqrt{n}$ }
\begin{center}
\begin{tabular}{|c|c|c|c|c|}
\hline
$c$ & 1/3 & 1/2 & 1 & 2 \\
\hline
\hline
n=$10^3$ & 1.4729  &  1.4176  &  1.3287 &   1.2458  \\
n=$10^4$ & 1.4352 &   1.3948 &   1.3287  &  1.2660 \\
n=$10^5$ &     1.4131  &  1.3813   & 1.3287  &  1.2783   \\
n=$10^6$& 1.3986   & 1.3723   & 1.3287  &  1.2866  \\
n=$10^7$ & 1.3884  &  1.3660 &   1.3287 &   1.2925 \\
n=$10^8$ & 1.3807 &   1.3613  &  1.3287  &  1.2970  \\
\hline
\hline
$c$ &  3 & 5 & 10  & 50\\
\hline
\hline
n=$10^3$   &  1.1997 &   1.1439  & 1.0716   & 0.9171\\
n=$10^4$ &   1.2307 &   1.1876  &  1.1313   & 1.0090\\
n=$10^5$   & 1.2497 &   1.2146 &   1.1685  &  1.0671\\
n=$10^6$ &  1.2626  &  1.2329 &   1.1938   & 1.1073\\
n=$10^7$  &   1.2718  &  1.2462   & 1.2122  &  1.1367\\
n=$10^8$  &  1.2788 &   1.2562   & 1.2262  &  1.1592\\
\hline
\end{tabular}
\end{center}
\label{OnSqrtTable}
\end{table}


\subsubsection{$\left|\mathcal{L}_n^0\right|=c\log n$}
 Again let $p_1=\mathbb{P}(N(\mu,1)>x^{\ast})$ and $\delta_2=0.0001$. Let $\mu$ be such that
\begin{equation}\label{LogMinMu}
\mu\sqrt{\log_{1/p_1}(c\log n)}>\sqrt{(2+\delta_2)\log mn}.
\end{equation}
We list the minimum value of $\mu$ that satisfies (\ref{LogMinMu}) in Table \ref{logNTable}.
In Table  \ref{logNTable}, we find that  the minimum detectable mean $\mu$ gradually increases as $n$ becomes larger.
This is due to the fact that the ratio $\left|\mathcal{L}_n^0\right|/n$ becomes more and more negligible as $n$ tends to $\infty$.
Table \ref{tb:ratioLog} gives the ratio of the length $\left|\mathcal{L}_n^0\right|$ of the embedded chain to the column number $n$ corresponding to the settings in Table  \ref{logNTable}.
When $n=10^8$ and $c=100$, the inhomogeneous chain only occupies about $1.8\times 10^{-5}$ portion of the images which is fairly negligible.
\begin{table}[htbp]
\caption{The minimum detectability of $\mu$ when $m=10$, $C=1$ and $\left|\mathcal{L}^0_n\right|=c\log n$ }
\begin{center}
\begin{tabular}{|c|c|c|c|c|c|c|}
\hline
$c$ & 1 & 2 & 5 & 10 & 50 & 100\\
\hline
n=$10^3$ & 1.83 & 1.70 & 1.58 & 1.51 & 1.39 & 1.35 \\
n=$10^4$ & 1.86 & 1.75 & 1.64 & 1.57 & 1.46 & 1.41\\
n=$10^5$ & 1.89 & 1.79 & 1.70 & 1.63 & 1.51 & 1.47\\
n=$10^6$ & 1.92 & 1.83 & 1.73 & 1.67 & 1.56 & 1.52\\
n=$10^7$ & 1.95 & 1.87 & 1.77 & 1.71 & 1.60 & 1.57\\
n=$10^8$ & 1.98 & 1.90 & 1.81 & 1.75 & 1.64 & 1.61\\
\hline
\end{tabular}
\end{center}
\label{logNTable}
\end{table}

\begin{table}[htbp]
\caption{The ratio of the length of the embedded chain to $n$.}
\begin{center}
\begin{tabular}{|c|c|c|c|}
\hline
$c$ & 1 & 2 & 5 \\
\hline
\hline
n=$10^3$ & 6.91$\times 10^{-3}$ & 1.38$\times 10^{-2}$ & 3.45$\times 10^{-2}$  \\
n=$10^4$ & 9.21$\times 10^{-4}$& 1.84$\times 10^{-3}$ & 4.61$\times 10^{-3}$ \\
n=$10^5$ & 1.15$\times 10^{-4}$ & 2.30$\times 10^{-4}$& 5.76$\times 10^{-4}$  \\
n=$10^6$ & 1.38$\times 10^{-5}$  & 2.76$\times 10^{-5}$ & 6.91e-$\times 10^{-5}$\\
n=$10^7$ & 1.16$\times 10^{-6}$ & 3.22$\times 10^{-6}$ & 8.06$\times 10^{-6}$\\
n=$10^8$ & 1.84$\times 10^{-7}$ & 3.68$\times 10^{-7}$ & 9.21$\times 10^{-7}$\\
\hline
\hline
$c$ & 10 & 50 & 100\\
\hline
\hline
n=$10^3$ & 6.91$\times 10^{-2}$ & 3.45$\times 10^{-1}$ & 6.91$\times 10^{-1}$ \\
n=$10^4$ & 9.21$\times 10^{-3}$ & 4.61$\times 10^{-2}$& 9.21$\times 10^{-2}$\\
n=$10^5$  & 1.15$\times 10^{-3}$  & 5.76$\times 10^{-3}$  & 1.15$\times 10^{-2}$ \\
n=$10^6$ & 1.38$\times 10^{-4}$& 6.91$\times 10^{-4}$ & 1.38$\times 10^{-3}$\\
n=$10^7$ & 1.61$\times 10^{-5}$ & 8.06$\times 10^{-5}$ & 1.61$\times 10^{-5}$\\
n=$10^8$ & 1.84$\times 10^{-6}$ & 9.21$\times 10^{-6}$ & 1.84$\times 10^{-5}$\\
\hline

\end{tabular}
\end{center}
\label{tb:ratioLog}
\end{table}

\subsection{Detection of solar flare}\label{subsec:solar}
In this subsection, we use a real data set based on the  solar data observatory to further verify our proposed testing procedures.
A solar flare is defined as a sudden, transient, and intense variation in brightness, which is usually observed over the Sun's surface \cite{augusto2011connection}.
The dataset is recorded in a video format and is publicly available online at \url{https://voices.uchicago.edu/willett/research/software/mousse/}.
There are in total 300 frames in the video, each of which contains a size of $232 \times 292$ image data.
According to the video, there are at least two obvious transient flares, which occur at frames $t = 187$$\sim$$202$ and $t = 216$$\sim$$268$,
respectively.
For illustrating purposes, the background information has been already removed and the remaining data is approximately normally distributed as mentioned in \cite{xie2013change}, where they used the first 100 frames to train the model.
In our experiment, we only
\begin{figure}[htbp]
\centering\includegraphics[width=100mm]{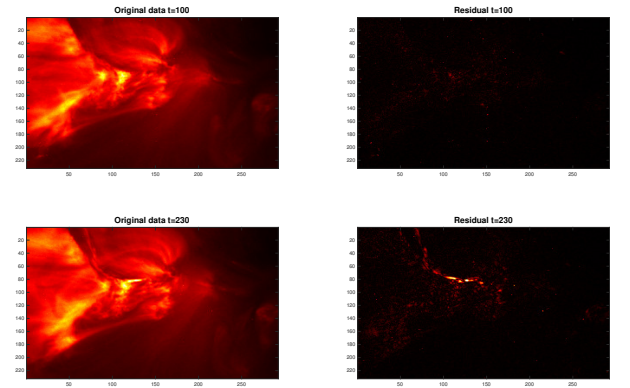}
\caption{Comparison of non-solar flare image and solar flare image: upper left shows a non-solar flare image; upper right shows a non-solar flare residual image; lower left shows a solar flare image; lower right shows a solar flare residual image.}
\label{fig_solar}
\end{figure}
compute the statistics from $t = 100$ and discard the frames used in the process of removing the background information.
In Fig. \ref{fig_solar}, we provide an example of the solar image taken at time $t = 100, 230$, respectively, along with the corresponding residual images. It can be clearly seen, there is a sudden intense burst in the middle of the residual image at $t = 230$ compared with no burst at all at $t = 100$.

\begin{figure}[htbp]
\centering\includegraphics[width=100mm]{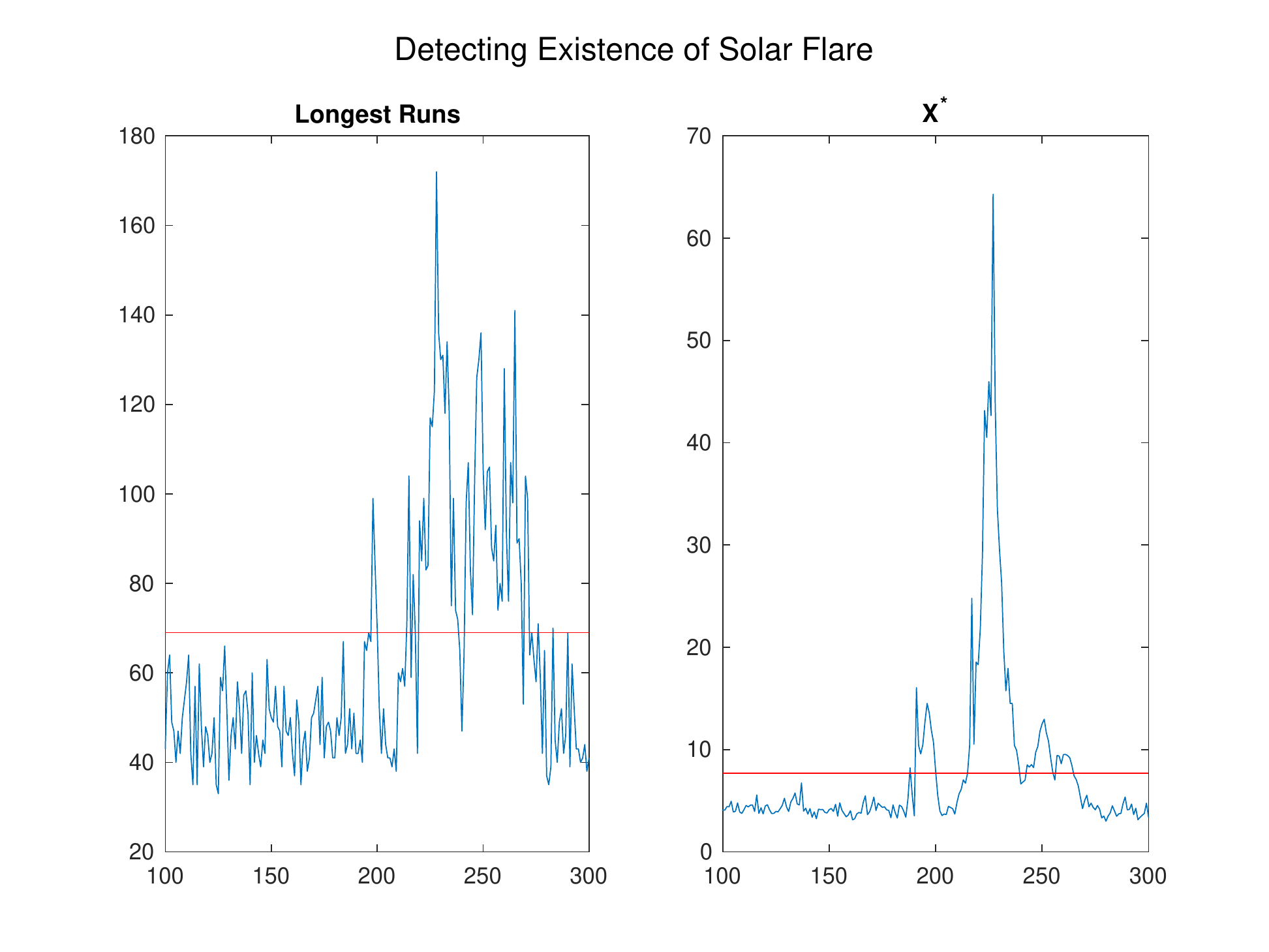}
\caption{Testing statistics from $t = 100$$\sim$$300$: left panel shows how our proposed testing statistics in Step 2 of the algorithm change; right panel shows how the length of the longest significant chain statistics in Step 1 of the algorithm change.}
\label{fig_solarStat}
\end{figure}

Naturally, the  corresponding burst region will have a potential significantly longer embedded chain than the images without a burst, like shown Fig. \ref{fig_solar}. And the value of the testing scan statistics $X^\ast$ at the occurrence of a solar flare will be larger than that of non-solar flare.
In Fig. \ref{fig_solarStat}, we compute the testing statistics of the length of the longest runs and the scan statistics for frames at $t = 100 \sim 300$.
It can be seen clearly that there are two big increase of the testing statistics around the time $t = 190$ and $t = 230$. By setting the threshold at $X^\ast_{thresh} = 7.7$ and $\left|L_0(n)\right|_{thresh} = 69$, our procedure can accurately identify the time of the occurrence of the solar flare, which corresponds to the type I error at level $0.05$.

%
%
%
%
%

%
%
%
%
%
%
%
%
%
%
%
\section{Extension}
\label{sec:Extension}
In this section, we will discuss about the longest significant run approach (as in Section \ref{sec:Summary}) in the detection problem of the $m$-by-$n$ array of nodes $\mathcal{S}$ as $m\rightarrow\infty$ and $n\rightarrow\infty$.
A model with similar structure is studied profoundly in \cite{ClustDetecPerco}.
After thresholding the values at the nodes with threshold $x^{\ast}$, under the null hypothesis $\mathbb{H}_0$, each node $(i,j)\in\mathcal{S}$ is significant with
\[
p=\mathbb{P}(N(0,1)>x^{\ast}).
\]
Let $p_1=\mathbb{P}(N(\mu,1)>x^{\ast})$ be the probability of significance under $\mathbb{H}_1(\mu)$. We use $L_0(m,n)$ to denote the longest chain consisting of significant nodes only under $\mathbb{H}_0$ and $\left|L_0(m,n)\right|$ is length.
In \cite{AssympConvergRate}, the authors show that there exists a continuous  function $\phi(p)$, which only depends on $C$ and $p$, but not on $m$ and $n$, such that as $m\rightarrow\infty$ and $n\rightarrow\infty$, we have
\begin{equation}\label{assympRateSub}
\frac{\left|L_0(m,n)\right|}{\log (mn)}\rightarrow\frac{1}{\phi(p)}, \text{\quad in probability},
\end{equation}
for $p<p_c$.
Here $p_c$ is a thresholding probability that the behaviors of the longest significant run are totally different when $p<p_c$ and $p\geq p_c$.
For more details, please refer to \cite{AssympConvergRate}.
Besides, $\phi(p)$ is a strictly decreasing function and it is positive when $p<p_c$ and constantly $0$ as $p\geq p_c$.
In \cite{AssympConvergRate}, it is shown that $p_c\geq\frac{1}{2C+1}$.
Thus, we may choose $x^{\ast}$ such that $p<\frac{1}{2C+1}$ under the null hypothesis.
As $m$ and $n$ become sufficiently large, the length of the longest significant chain is at most $(1+\epsilon)\frac{\log mn}{\phi(p)}$ for some $\epsilon>0$.
Given the above, we have the following revised detection algorithm for the case that $(m,n)\rightarrow(\infty,\infty)$.\\
\textbf{Detection Algorithms when $(m,n)\rightarrow(\infty,\infty)$:}
\begin{enumerate}
\item Take $x^{\ast}$ such that $p=\mathbb{P}(N(0,1)>x^{\ast})<\frac{1}{2C+1}$ to be the threshold of nodes to be significant.
Let $\mathcal{E}_n=\{\mathcal{L}\in\mathcal{F}_n: Z(\mathcal{L})=1\}$.
Find the longest chain $L_0(m,n)$ in $\mathcal{E}_n$.
For small $\epsilon>0$ if the length $|L_0(m,n)|>(1+\epsilon/2)\frac{\log(mn)}{\phi(p)}$, then reject $\mathbb{H}_0$; otherwise, go to the next step.
\item Compute $X^{\ast}_s$ as
\[
\max_{\mathcal{L}\in\mathcal{E}_n}\sum_{(i,j)\in\mathcal{L}}\frac{X(i,j)}{\sqrt{\left|\mathcal{L}\right|}}.
\]
For small $\delta_2>0$, if $X^{\ast}_s>\sqrt{2(1+\delta_2)\log(mn)}$, then reject $\mathbb{H}_0$; otherwise accept $\mathbb{H}_0$.
\end{enumerate}

As shown above, the first step takes $O(mn)$ and second $O(mn\log(mn))$.
Hence this algorithm takes $O(mn\log(mn))$ flops in total with $O(mn\log(mn))$ required space for storage.
Moreover, by Theorem \ref{thLongSigChain} and \ref{type2error}, it is straightforward to see that under $\mathbb{H}_1$,
\begin{enumerate}
\item if $\left|\mathcal{L}^0_n\right|\geq\zeta_1 n$ for some $\zeta_1>0$ and $p_1>\exp\{-\phi(p)\frac{\log\zeta_1 n}{(1+\epsilon)\log(mn)}\}$;
\item or if $\mu\sqrt{\log_{1/p_1}\left|\mathcal{L}^0_n\right|}>\sqrt{2(1+\delta_2)\log(mn)}$ for some $\delta_2$,
\end{enumerate}
then, as $m\rightarrow\infty$ and $n\rightarrow\infty$, we have
\begin{equation*}
\mathbb{P}(\text{accept } \mathbb{H}_1\big|\mathbb{H}_0)+\mathbb{P}(\text{accept } \mathbb{H}_0\big|\mathbb{H}_1)\rightarrow 0.
\end{equation*}

\section{Conclusion}
In this paper, we give a detection method for chains with elevated means in a white noise image.
We analyze the length of the longest significant chain after thresholding each pixel and consider the statistics over all significant chains.
Such a strategy significantly reduces the complexity of the algorithm and the false positives are eliminated as the number of pixels increases.
The numeric study shows the results are very promising, compared to human eyes' detectability.
The real data example on solar flare detection also verifies the effectiveness of our proposed method.

\appendices

\section*{Appendix}
This Appendix includes the proofs of our main technical results in Section \ref{sec:LongSigRun}.
Specifically, in Appendix \ref{app:a}, we provide the proof for the Corollary \ref{co:depofrho}, which states the relation between $\rho$ and $(m, p)$.
The proof of Theorem \ref{type1error}, which shows under $\mathbb{H}_0$, that the type-I error tends to $0$ as $n\rightarrow\infty$, is in Appendix \ref{app:b}.
The proofs for Theorem \ref{thLongSigChain} and \ref{type2error} regarding the asymptotic diminishing rate of the type-II error, are provided in Appendix \ref{app:c} and \ref{app:d}, respectively.

\section{Proof of Corollary \ref{co:depofrho}}\label{app:a}
Given a realization
$$
t_{i,j}\sim N(0,1), 1\leq i\leq n, 1\leq j\leq m,
$$
let $x^{\ast}_1\geq x^{\ast}_2>0$ be such that $p_1=\mathbb{P}(t_{i,j}>x^{\ast}_1)$ and $p_2=\mathbb{P}(t_{i,j}>x^{\ast}_2)$.
Since $t_{i,j}>x^{\ast}_1$ implies that $t_{i,j}>x^{\ast}_2$, one can easily see that each significant node under threshold $x^{\ast}_1$ must be significant under $x^{\ast}_2$.
  Therefore, we have $\left|L_0(m_1,n,p_1)\right|\leq\left|L_0(m_1,n,p_2)\right|$, which by Theorem \ref{thConvRate} leads to $\rho(m_1,p_1)\leq \rho(m_1,p_2)$.
Similarly, it is not hard to see that $\left|L_0(m_1,n,p_1)\right|\leq\left|L_0(m_2,n,p_1)\right|$, since if $m_1\leq m_2$, $([1,n]\times[1,m_1])\cap\mathbb{Z}^2\subset([1,n]\times[1,m_2])\cap\mathbb{Z}^2$.
Thus $\rho(m_1,p_1)\leq\rho(m_2,p_1)$.

\section{Proof of Theorem \ref{type1error}}\label{app:b}
Recall (\ref{eqEgoroff}) that for any $\delta>0$ and $\epsilon>0$, the length of the longest significant chain under $\mathbb{H}_0$ is no larger than $b_{\epsilon, n}=(1+\epsilon)\log_{1/\rho}n$ with probability $1-\delta$.
By the Bonferroni inequality, it is easy to derive the following,

\begin{eqnarray}
&&\mathbb{P}(X_s^{\ast}>\tau)\nonumber\\
&=&\mathbb{P}(\bigcup_{\mathcal{L}\in\mathcal{F}_n}\{X(\mathcal{L})>\tau,Z(\mathcal{L})=1\})\nonumber\\
&\leq&\delta+\mathbb{P}(\bigcup_{k=1}^{b_{\epsilon,n}}\bigcup_{\mathcal{L}\in\mathcal{F}_n, \left|\mathcal{L}\right|=k}\{X(\mathcal{L})>\tau, Z(\mathcal{L})=1\})\nonumber\\
&\leq&\delta+\sum_{k=1}^{b_{\epsilon,n}}\sum_{\mathcal{L}\in\mathcal{F}_n, \left|\mathcal{L}\right|=k}\mathbb{P}(X(\mathcal{L})>\tau\big|Z(\mathcal{L})=1)\times \nonumber \\&&\mathbb{P}(Z(\mathcal{L})=1)\nonumber\\
&\leq&\delta+\sum_{k=1}^{b_{\epsilon,n}}\sum_{\mathcal{L}\in\mathcal{F}_n, \left|\mathcal{L}\right|=k}p^k\mathbb{P}(X(\mathcal{L})>\tau\big|Z(\mathcal{L})=1),
\label{continueEq}
\end{eqnarray}

where $\left|\mathcal{L}\right|$ is the length of the significant chain $\mathcal{L}$.
Note that conditioning on the event $Z(\mathcal{L})=1$, each $X(i,j)$ on $\mathcal{L}$  is a truncated standard normal random variable bounded below by $x^{\ast}$.
Let $X_T(i,j)$ ($1\leq i\leq n, 1\leq j\leq m$) be i.i.d. random variables with distribution equal to $X(i,j)$ given that $X(i,j)>x^{\ast}$, where $X(i,j)$ is the standard normal random variable under $\mathbb{H}_0$.
Let $\Phi(\cdot)$ be the distribution function of the standard normal distribution.
It is easy to see that the probability density function for $X_T(i,j)$ is $f_{X_T}(x) = \frac{1}{(1-\Phi(x^{\ast}))\sqrt{2\pi}}\exp\{-\frac{x^2}{2}\}$.
Thus, for each $\mathcal{L}\in\mathcal{E}_n$, we have
\begin{eqnarray}
&&\mathbb{P}(X(\mathcal{L})>\tau\big|Z(\mathcal{L})=1)\nonumber\\
&=&\mathbb{P}(\sum_{(i,j)\in\mathcal{L}}\frac{X(i,j)}{\sqrt{\left|\mathcal{L}\right|}}>\tau\big|Z(\mathcal{L})=1)\nonumber\\
&=&\mathbb{P}(\sum_{(i,j)\in\mathcal{L}}X_{T}(i,j)>\tau\sqrt{\left|\mathcal{L}\right|})\nonumber\\
&\displaystyle \stackrel{\text{(Markov's Inequality)}}\leq&\inf_{\omega\geq0}\exp\{-\omega\tau\sqrt{\left|\mathcal{L}\right|}\}\prod_{(i,j)\in\mathcal{L}}\mathbb{E}\exp\{\omega X_{T}(i,j)\}\nonumber\\
&\displaystyle \stackrel{\text{(p.d.f. of $X_{T}(i,j)$)}}\leq&\inf_{\omega\geq0}\exp\{-\omega\tau\sqrt{\left|\mathcal{L}\right|}\}\nonumber \times \prod_{(i,j)\in\mathcal{L}}\frac{1}{(1-\Phi(x^{\ast}))\sqrt{2\pi}}\nonumber\\
&&\int_{x^{\ast}}^{\infty}\exp\{\omega y-\frac{y^2}{2}\}\mathrm{d}y\nonumber\\
&=&\inf_{\omega\geq0}\exp\{-\omega\tau\sqrt{\left|\mathcal{L}\right|}\}\nonumber
\times \prod_{(i,j)\in\mathcal{L}}\frac{\exp\{\frac{\omega^2}{2}\}}{(1-\Phi(x^{\ast}))\sqrt{2\pi}}\nonumber\\
&&\int_{x^{\ast}}^{\infty}\exp\{-\frac{(y-\omega)^2}{2}\}\mathrm{d}y\nonumber\\
&=&\inf_{\omega\geq0}\exp\{-\omega\tau\sqrt{\left|\mathcal{L}\right|}\}\nonumber \times \prod_{(i,j)\in\mathcal{L}}\frac{\exp\{\frac{\omega^2}{2}\}}{(1-\Phi(x^{\ast}))\sqrt{2\pi}}\nonumber\\
&&\int_{x^{\ast}+\omega}^{\infty}\exp\{-\frac{z^2}{2}\}\mathrm{d}z\nonumber\\
&=&\inf_{\omega\geq0}\exp\{-\omega\tau\sqrt{\left|\mathcal{L}\right|}\}\nonumber \times \prod_{(i,j)\in\mathcal{L}}\frac{\exp\{\frac{\omega^2}{2}\}}{(1-\Phi(x^{\ast}))}\nonumber\\
&&(1-\Phi(x^{\ast}+\omega))\nonumber\\
&\leq&\inf_{\omega\geq0}\exp\{-\omega\tau\sqrt{\left|\mathcal{L}\right|}+\left|\mathcal{L}\right|\frac{\omega^2}{2}\}\nonumber\\
&=&\inf_{\omega\geq0}\exp\{\frac{1}{2}(\omega\sqrt{\left|\mathcal{L}\right|}-\tau)^2-\frac{(\tau)^2}{2}\}\nonumber\\
&\leq&\exp\{-\frac{(\tau)^2}{2}\}.\label{pluginContEq}
\end{eqnarray}
Plug (\ref{pluginContEq}) into (\ref{continueEq}), since by our assumption that $p=\mathbb{P}(N(0,1)>x^{\ast})\leq\rho<\frac{1}{2C+1}$, we have
\begin{eqnarray*}
&&\mathbb{P}(X_{s}^{\ast}>\tau)\\
&\leq&\delta+\sum_{k=1}^{b_{\epsilon,n}}\sum_{\mathcal{L}\in\mathcal{F}_n, \left|\mathcal{L}\right|=k}p^k\exp\{-\frac{(\tau)^2}{2}\}\\
&\leq&\delta+mn\exp\{-\frac{(\tau)^2}{2}\}\sum_{k=1}^{b_{\epsilon,n}}(2C+1)^{k}p^{k}\\
&\leq&\delta+mn\frac{1}{1-(2C+1)p}\exp\{-\frac{(\tau)^2}{2}\}.
\end{eqnarray*}
Since $\delta>0$ is arbitrary and $m$ is fixed, if $\tau_s^{\ast}=\sqrt{2(1+\delta_2)\log mn}$, then we have that $\mathbb{P}(X_s^{\ast}>\tau_s^{\ast}\big|\mathbb{H}_0)\rightarrow0$ for any small $\delta_2>0$, as $n\rightarrow\infty$. 

\section{Proof of Theorem \ref{thLongSigChain}}\label{app:c}
By the argument before the theorem, for any small $\delta>0$ and $\epsilon_1>0$, with probability $1-\delta$, there exists $N\in\mathbb{Z}^{+}$ such that when $n\geq N$,
\[
\left|L_0(n)\right|\geq\left|L_1(n)\right|\geq(1-\epsilon_1)\log_{1/p_1}\left|\mathcal{L}^0_n\right|.
\]
Under $\mathbb{H}_1$, by (\ref{condOfLongSigRun}), one can derive
\[
\log_{1/p_1}\left|\mathcal{L}_n^0\right|\geq \alpha\log_{1/p_1}\zeta_1 n>(1+\epsilon)\log_{1/\rho}n.
\]
We choose $\epsilon_1$ such that $(1+\epsilon)(1-\epsilon_1)>1+\epsilon/2$ and thus
\[
\mathbb{P}[\left|L_0(n)\right|>(1+\epsilon/2)\log_{1/\rho}n\big|\mathbb{H}_1]\geq 1-\delta, \forall n\geq N.
\]
Since $\delta$ is arbitrary, we have (\ref{asympPowerOfLongSigRun}) asymptotically.
The second part of (\ref{condOfLongSigRun}) follows from the fact that
\[
\lim_{n\to\infty}\frac{\log\zeta n}{(1+\epsilon)\log n}=\lim_{n\to\infty}\frac{\log\zeta+\log n}{(1+\epsilon)\log n}=\frac{1}{1+\epsilon}.
\]

\section{Proof of Theorem \ref{type2error}}\label{app:d}
Before the proof, let us first  recall the following definition about the association of random variables in \cite{AORV}.
\begin{definition}\label{association}
We say random variables $T_1, T_2, \ldots, T_n$ are associated, if $\text{Cov}[f(\bf{T}), g(\bf{T})]\geq 0,$ where $\bf{T}=(T_1, T_2, \ldots, T_n)$, for any nondecreasing functions $f$ and $g$, for which $\mathbb{E}f(\bf{T})$, $\mathbb{E}g(\bf{T})$, and $\mathbb{E}f(\bf{T})g(\bf{T})$ exist.
\end{definition}

Now we give the proof of Theorem \ref{type2error} in the following.

Recall in (\ref{egoroffunderH1}), with high probability, the length of  the longest significant chain in $\mathcal{L}^0_n$, $\left|L_1(n)\right|$ falls in the region: $[c^{L}_{\epsilon,n},c^{U}_{\epsilon,n}]$. It is not difficult to derive the following,
\begin{eqnarray}
&&\mathbb{P}(X_s^{\ast}>\tau_s^{\ast}|\mathbb{H}_1)\nonumber\\
&=&\mathbb{P}(\bigcup_{\mathcal{L}\in\mathcal{F}_n}\{X(\mathcal{L})>\tau_s^{\ast}, Z(\mathcal{L})=1\}|\mathbb{H}_1)\nonumber\\
&\geq&\mathbb{P}(X(L_1(n))>\tau_s^{\ast})\nonumber\\
&=&\sum_{k=1}^{n}\mathbb{P}(\sum_{(i,j)\in L_1(n)}\frac{X(i,j)}{\sqrt{\left|L_1(n)\right|}}>\tau_s^{\ast}\big|\left|L_1(n)\right|=k)\nonumber \times\nonumber\\
&&\mathbb{P}(\left|L_1(n)\right|=k)\nonumber\\
&\geq&\sum_{k=c^{L}_{\epsilon,n}}^{c^{U}_{\epsilon,n}}\mathbb{P}(\sum_{(i,j)\in L_1(n)}\frac{X(i,j)}{\sqrt{k}}>\tau_s^{\ast}\big| \left|L_1(n)\right|=k)\nonumber \times \nonumber\\
&&\mathbb{P}(\left|L_1(n)\right|=k)\label{ineq:split}
\end{eqnarray}

Recall that $\mathcal{F}_n$ is the set of all chains of good continuation in $\mathcal{S}=\{(i,j): 1\leq i\leq n, 1\leq j\leq m\}.$
Let $\mathcal{F}^k_n\subset\mathcal{F}_n$ be the set of all chains of good continuation with length $k$, namely, $\mathcal{F}^k_n=\{\mathcal{L}\in\mathcal{F}_n: \left|\mathcal{L}\right|=k\}.$

One can easily see that
\begin{eqnarray*}
&&\mathbb{P}(\sum_{(i,j)\in L_1(n)}\frac{X(i,j)}{\sqrt{k}}>\tau_s^{\ast}\big| \left|L_1(n)\right|=k)\\
&=&
\frac{\sum_{\mathcal{L}\in\mathcal{F}^k_n}\mathbb{P}(\sum_{(i,j)\in\mathcal{L}}\frac{X(i,j)}{\sqrt{k}}>\tau_s^{\ast}\big| L_1(n)=\mathcal{L})\mathbb{P}(L_1(n)=\mathcal{L})}{\sum_{\mathcal{L}\in\mathcal{F}^k_n}\mathcal{P}(L_1(n)=\mathcal{L})}.
\end{eqnarray*}

Generate $k$ random variables $Y_1, \ldots, Y_k\stackrel{\text{i.i.d.}}{\sim}N(\mu, 1)$ which are independent from $\{X_{i,j}, 1\leq i\leq n, 1\leq j\leq m\}$. For each $\mathcal{L}\in\mathcal{F}^k_n$, we have that
\begin{eqnarray*}
&&\mathbb{P}(\sum_{(i,j)\in\mathcal{L}}\frac{X(i,j)}{\sqrt{k}}>\tau^{\ast}_s\big| L_1(n)=\mathcal{L})\\
&=&\mathbb{P}(\sum_{i=1}^{k}Y_i\big/\sqrt{k}>\tau^{\ast}_s\big| Y_1>x^{\ast},\ldots, Y_k>x^{\ast}).
\end{eqnarray*}

Let $A$ and $B$ be two subsets of $\mathbb{R}^{k}$ such that
\[
A=\{(y_1,\ldots,y_k): \sum_{i=1}^{k}\frac{y_i}{\sqrt{k}}>\tau^{\ast}_s\},
\]
\[
B=\{(y_1,\ldots,y_k):y_1>x^{\ast},\ldots,y_k>x^{\ast}\}.
\]
Let $f:\mathbb{R}^k\rightarrow\{0,1\}$ and $g:\mathbb{R}^k\rightarrow\{0,1\}$ be indicator functions of sets $A$ and $B$, respectively, i.e.,
\[
f(y_1,\ldots,y_k)=I_{A}(y_1,\ldots,y_k)
\]
and
\[
g(y_1,\ldots,y_k)=I_{B}(y_1,\ldots,y_k).
\]

Theorem 2.1 of \cite{AORV} states that independent random variables are associated.
Therefore $Y_1, \ldots, Y_k$ are associated since they are independent.
Realize that both $f$ and $g$ are increasing functions and therefore, it is straightforward to see that
\begin{eqnarray*}
&&\mathbb{P}(\sum_{i=1}^k Y_i\big/\sqrt{k}>\tau^{\ast}_s, Y_1>x^{\ast},\ldots,Y_k>x^{\ast})\\
&=&\mathbb{E}[f(Y_1,\ldots,Y_k)g(Y_1,\ldots,Y_k)]\\
&\geq&\mathbb{E}f(Y_1,\ldots,Y_k)\mathbb{E}g(Y_1,\ldots,Y_k)\\
&=&\mathbb{P}(\sum_{i=1}^k Y_i/\sqrt{k}>x^{\ast})\mathbb{P}(Y_1>x^{\ast},\ldots,Y_k>x^{\ast}).
\end{eqnarray*}
Hence it follows that
\begin{eqnarray*}
&&\mathbb{P}(\sum_{i=1}^k Y_i\big/\sqrt{k}>\tau^{\ast}_s\big| Y_1>x^{\ast},\ldots,Y_k>x^{\ast})\\
&\geq&\mathbb{P}(\sum_{i=1}^k Y_i\big/\sqrt{k}>\tau^{\ast}_s)\\
&=&\mathbb{P}(N(\sqrt{k}\mu,1)>\tau_s^{\ast}).
\end{eqnarray*}
Since $k\geq c^{L}(\epsilon,n)$, as long as
\[
\mu\sqrt{c^L(\epsilon, n)}>\tau_s^{\ast}=\sqrt{2(1+\delta_2)\log n},
\]
by Mill's ratio, we have
\begin{eqnarray*}
&&\mathbb{P}(N(\sqrt{k}\mu,1)<\tau_s^{\ast})\\
&\leq&\mathbb{P}(N(0,1)<-\gamma \sqrt{\log n})\\
&\leq&2n^{-\gamma^2/2}\rightarrow0, \text{ as } n\rightarrow 0,
\end{eqnarray*}
where
\begin{eqnarray*}
\gamma&=&\sqrt{k}\mu\big/\sqrt{\log n}-\sqrt{2(1+\delta_2)}\\
&\geq&\sqrt{c^L(\epsilon, n)}\mu\big/\sqrt{\log n}-\sqrt{2(1+\delta_2)} > 0.
\end{eqnarray*}
Therefore, going back to (\ref{ineq:split}), it follows that
\begin{eqnarray*}
\mathbb{P}(X_s^{\ast}>\tau_s^{\ast})&\geq&\sum_{k=c^{L}_{\epsilon,n}}^{c^U_{\epsilon,n}}(1-2n^{-\gamma^2/2})\mathbb{P}(\left|L_1(n)\right|=k)\\
&\geq&(1-2n^{-\gamma^2/2})(1-\delta).
\end{eqnarray*}
Since $\delta>0$ is arbitrary, as $n\rightarrow\infty$, we have
\[
\mathbb{P}(X^{\ast}_s>\tau_s^{\ast}\big|\mathbb{H}_1)\rightarrow 1.
\]


\section*{Acknowledgment}
The authors would like to thank Professor Vladimir Koltchinskii for his careful reading on this manuscript.
The first author would like to express his gratitude to Professor Koltchinskii for his continuous support on his Ph.D. study in Georgia Institute of Technology.
This project is partially supported by the Transdisciplinary Research Institute for Advancing Data Science (TRIAD), http://triad.gatech.edu, which is a part of the TRIPODS program at NSF and locates at Georgia Tech, enabled by the NSF grant CCF-1740776.
Cao and Huo are partially supported by the NSF grant 1613152.

\ifCLASSOPTIONcaptionsoff
  \newpage
\fi



\bibliographystyle{IEEEtran}
\bibliography{lrrefnew}

%

%
%
%




\end{document}